\documentclass[12pt]{amsart}
\usepackage{amssymb, amstext, amscd, amsmath}
\makeatletter
\def\@cite#1#2{{\m@th\upshape\bfseries%
[{#1\if@tempswa{\m@th\upshape\mdseries, #2}\fi}]}}
\makeatother
\theoremstyle{plain}
\newtheorem{thm}{Theorem}[section]
\newtheorem{cor}[thm]{Corollary}
\newtheorem{prop}[thm]{Proposition}
\newtheorem{lem}[thm]{Lemma}
\theoremstyle{definition}
\newtheorem{rem}[thm]{Remark}
\newtheorem{defn}[thm]{Definition}

\newtheorem{eg}[thm]{Example}
\newcommand{\Prf}{\noindent\textbf{Proof.\ }}
\newcommand{\bx}{\hfill$\blacksquare$\medbreak}

\newcommand{\bsl}{\setminus}

\newcommand{\lip}{\langle}
\newcommand{\rip}{\rangle}
\newcommand{\ip}[1]{\lip #1 \rip}
\newcommand{\mt}{\emptyset}
\newcommand{\ol}{\overline}

\DeclareMathOperator*{\sotlim}{\textsc{sot}--lim}
\newcommand{\sotsum}{\textsc{sot--}\!\!\!\!\sum}
\newcommand{\wot}{\textsc{wot}}

\newcommand{\bC}{{\mathbb{C}}}
\newcommand{\bF}{{\mathbb{F}}}

\newcommand{\bN}{{\mathbb{N}}}

\newcommand{\bQ}{{\mathbb{Q}}}
\newcommand{\bR}{{\mathbb{R}}}
\newcommand{\bT}{{\mathbb{T}}}
\newcommand{\bZ}{{\mathbb{Z}}}

  \newcommand{\B}{{\mathcal{B}}}
  \newcommand{\C}{{\mathcal{C}}}
  \newcommand{\D}{{\mathcal{D}}}
  \newcommand{\E}{{\mathcal{E}}}

\renewcommand{\H}{{\mathcal{H}}}
  
  \newcommand{\J}{{\mathcal{J}}}
  \newcommand{\K}{{\mathcal{K}}}
\renewcommand{\L}{{\mathcal{L}}}
  \newcommand{\M}{{\mathcal{M}}}
  \newcommand{\N}{{\mathcal{N}}}
\renewcommand{\O}{{\mathcal{O}}}
\renewcommand{\P}{{\mathcal{P}}}
  
  \newcommand{\R}{{\mathcal{R}}}

  \newcommand{\T}{{\mathcal{T}}}
  \newcommand{\U}{{\mathcal{U}}}
  \newcommand{\V}{{\mathcal{V}}}

\newcommand{\upchi}{{\raise.35ex\hbox{$\chi$}}}

\newcommand{\fA}{{\mathfrak{A}}}

\newcommand{\fK}{{\mathfrak{K}}}
\newcommand{\fL}{{\mathfrak{L}}}

\newcommand{\fR}{{\mathfrak{R}}}
\newcommand{\fS}{{\mathfrak{S}}}

\newcommand{\AND}{\text{ and }}
\newcommand{\FOR}{\text{ for }}

\newcommand{\qand}{\quad\text{and}\quad}
\newcommand{\qandal}{\quad\text{and all}\quad}
\newcommand{\qif}{\quad\text{if}\quad}
\newcommand{\qfor}{\quad\text{for}\quad}
\newcommand{\qforal}{\quad\text{for all}\quad}

\newcommand{\Alg}{\operatorname{Alg}}

\newcommand{\diag}{\operatorname{diag}}

\newcommand{\ran}{\operatorname{Ran}}
\newcommand{\rank}{\operatorname{rank}}
\newcommand{\spn}{\operatorname{span}}
\newcommand{\sumoplus}{\sum^\oplus}

\newcommand{\fgeeplus}{\bF^+(G)}

\newcommand{\rowt}{(T_e)_{e\in \E(G)}}

\newcommand{\flgee}{\fL_G}

\begin{document}

\title[Nest representations]{Nest representations\\of directed graph algebras}
%
\author[K.R.Davidson]{Kenneth R. Davidson}
\thanks{First author partially supported by an NSERC grant.}
\address{Pure Math.\ Dept.\\U. Waterloo\\Waterloo, ON\;
N2L--3G1\\CANADA}
\email{krdavids@uwaterloo.ca}
\author[E. Katsoulis]{Elias~Katsoulis}
\address{Dept. Math.\\East Carolina University\\
Greenville, NC 27858\\USA}
\email{KatsoulisE@mail.ecu.edu}
\begin{abstract}

This paper is a comprehensive study of the nest
representations for the free semigroupoid algebra $\flgee$
of countable directed graph $G$ as well as its norm-closed
counterpart, the tensor algebra $\T^{+}(G)$.

We prove that the finite dimensional nest representations separate the
points in $\flgee$, and a fortiori, in $\T^{+}(G)$.
The irreducible finite dimensional representations separate the
points in $\flgee$ if and only if $G$ is transitive in components (which
is equivalent to being semisimple).
Also the upper triangular nest representations separate points if and
only if for every vertex $x \in \V(G)$ supporting a cycle, $x$ also
supports at least one loop edge.

We also study \textit{faithful} nest representations. We prove
that $\flgee$ (or $\T^{+}(G)$) admits a faithful irreducible
representation if and only if $G$ is strongly transitive as a directed
graph.  More generally, we obtain a condition on $G$ which is equivalent
to the existence of a faithful nest representation. We also give a
condition that determines the existence a faithful nest representation
for a maximal type $\bN$ nest.
\end{abstract}

%
\thanks{June 16, 2005 version.}
\date{}
\maketitle

\section{Introduction}

The representation theory for non-selfadjoint operator algebras,
although still in its infancy, is already playing an important role.
Given a non-selfadjoint algebra,  its nest representations are thought of
as the proper generalization of the irreducible ones and a general theory
around them is developing \cite{Lam, DKPet, KaPet}.
This paper contributes to the further development of that theory by
analyzing such representations for algebras associated to directed
graphs.

The study of operator algebras associated to directed graphs is fairly
recent in both C*-algebras and non-selfadjoint algebras.
While we do not use the C*-algebraic theory here, we refer to a few papers for
the interested reader \cite{KPR,KPRR,BPRS,DT}.
The non-selfadjoint theory is even more recent.  The finite graphs yield
quiver algebras, which have been studied by algebraists for some time and
were introduced to operator algebras by Muhly \cite{M}.
Much of the motivations for questions in this area have come from the
case of free semigroup algebras, which corresponds to a graph consisting
of a single vertex and $n$ loop edges \cite{DP1,DP2,DP3,DKS,DKP,DLP}.
Non-selfadjoint algebras of infinite graphs have been studied by a number of
authors beginning with Kribs and Power \cite{KrP}.
See also \cite{KK,JK1,JK2,KrP2}.

In \cite{KK}, the quiver and free semigroupoid algebras were
classified by making an essential use of their two dimensional nest
representations (see also \cite{Sol}). The analysis there raises the question: to what
extent do the finite dimensional nest representations determine the
algebra. Questions of this kind have been investigated for other
classes of operator algebras but apparently little has been done here,
even for the free semigroup algebras $\fL_n$, $n \geq2$.

In this paper, we prove that for an arbitrary countable directed graph $G$,
the finite dimensional nest representations of $\T^{+}(G)$ or
$\flgee$ separate points (Theorem \ref{fdnest_separate}).
If $G$ is transitive in components, then Theorem \ref{irreducible} shows that the
\textit{irreducible} finite dimensional representations alone are
sufficient to separate the
points.  Indeed this property characterizes transitivity
in components for $G$. At the other extreme, Theorem \ref{upper} shows that the upper
triangular nest representations separate the points if and only if for
every vertex $x \in \V(G)$ supporting a cycle, $x$ also supports at
least one loop edge. Both these results apply to the non-commutative disc algebra
$\fA_n$ and the non-commutative Toeplitz algebra $\fL_n$, for $n
\geq2$, thus showing that these algebras are strongly semisimple, i.e. the
maximal ideals have zero intersection.

Recently, Read \cite{Read} proved that $B(\H)$ is a free semigroup
algebra generated by two isometries with orthogonal ranges.
(See also \cite{D_read}.) In the language of representation theory, his result
asserts that the non-commutative disc algebra $\fA_2$ admits a faithful
$*$-extendible irreducible representation. If one drops the $*$-extendibility,
things simplify considerably.  We show that $\fL_G$ (or $\T^{+}(G)$)
admits a faithful irreducible representation if and only if the graph
$G$ is strongly transitive. This result is new even for the free
semigroup algebras $\fL_n$, $n\geq 2$.

The irreducible representations are a special case of nest representations.
The above considerations suggest the problem of characterizing which
graph algebras admit faithful nest representations. It turns out that not all
graph algebras admit faithful nest representations. Indeed, in Theorem \ref{faithful_nest_rep}
we identify an intricate set of conditions on
the graph which determine the existence of a faithful nest representation. In particular,
the transitive quotient of such a graph $G$ should inherit a total ordering from the order of $G$.
We also characterize (Theorem \ref{Nrepn}) when there exists a faithful nest representation
on the upper triangular operators, with respect to a basis
ordered as $\bN$.  If $G$ has no sinks, this characterization simply requires that the graph is transitive
and each one of its edges supports at least one loop edge. In particular, the free semigroup algebras
$\fL_n$, $n \ge 2$, admit such representations.

\begin{rem} \label{R:tensor}
In what follows, we state and prove the above mentioned results only
in the case of a free semigroupoid algebra. All these results are
valid in the tensor algebra case with identical statements and
proofs. See Proposition~\ref{P:tensor}.
We therefore omit them for brevity of the presentation.
\end{rem}

\section{Notation and Background}

A countable directed graph $G$ is given by $(\V(G), \E(G), r, s)$ where
$\V(G)$ is a countable set of objects called vertices,
$\E(G)$ is a countable set of objects called edges,
and $r$ and $s$ are maps from $\E(G)$ into $\V(G)$.
The maps are called the range and source maps, and we think of
$G$ as a set of points $x \in \V(G)$ with directed edges $e \in \E(G)$
connecting the vertex $s(e)$ to $r(e)$.
A vertex $x$ is called a sink of there is no edge $e$ with $s(e)=x$,
and a source if there is no edge $e$ with $r(e) = x$.

A path in $G$ is a finite sequence of edges which connect, range to
source, within the graph,
that is, a finite sequence $p = e_k e_{k-1} \dots e_2 e_1$ such that
$r(e_i) = s(e_{i+1})$ for $1 \le i < k$.
We will write these paths from right to left because this will
coincide with operator products in the representations of the graphs.
The length of a path is the number of edges, and is written $|p|$.
A vertex $x$ will be considered to be a path of zero length with source and
range maps $s(x) = r(x) = x$.

If $p = e_k e_{k-1} \dots e_2 e_1$ is a path, we set $r(p) = r(e_k)$ and
$s(p) = s(e_1)$. A cycle is a path $p$ such that $r(p) = s(p)$.
A loop edge is an edge $e$ which is a cycle, that is $r(e) = s(e)$.

A graph $G$ is \textit{transitive} if there is a path from each vertex in
the graph to every other.
Say that a graph $G$ is \textit{strongly transitive} if it is transitive
and is not a single cycle nor is trivial (no edges, one vertex).
The transpose of $G$, $G^t$, has the same same vertices and edges
as $G$, but the source and range maps are reversed: $r_{G^t} = s_G$
and $s_{G^t} = r_G$.
The graph $G$ is connected if the undirected graph determined by $G$
is transitive.
A connected component of $G$ is a maximal connected subset of $G$.
We say that $G$ is transitive on components is each component of $G$ is
transitive.
We will also talk of the \textit{transitive component} of a vertex $x$ to mean the
largest transitive subgraph of $G$ containing $x$.

The path space $\fgeeplus$ consists of all finite paths including the vertices.
We may consider $\fgeeplus$ as a semigroupoid with a product $pq$
defined if $r(q) = s(p)$.
In particular, $p = p s(p) = r(p) p$.

A representation of $G$ is a norm continuous homomorphism $\rho$ from
$\fgeeplus$ into $\B(\H)$ for some Hilbert space $\H$.
Since the vertices are idempotents in $\fgeeplus$ which annihilate
each other, a well-known result of Dixmier implies that,
after a similarity, we may assume that $\rho(x)$ are pairwise
orthogonal \textit{projections} for $x \in \V(G)$.
The representation is non-degenerate if $\H$ is spanned by the
ranges of these projections.
If $p$ is an edge, then it is easy to verify that
$\rho(p) = \rho(r(p)) \rho(p) \rho(s(p))$.
It turns out that the representations constructed
in this paper are
\textit{completely contractive}:
this is equivalent to the row operator
\[
 \begin{bmatrix}\rho(e_1) & \rho(e_2) & \dots\end{bmatrix}_{e_i\in\E(G)}
\]
being a contraction.  Of course one could demand instead that the matrices
$\begin{bmatrix}\rho(e_1) & \dots & \rho(e_k)\end{bmatrix}$
be contractions for every finite subset of distinct edges of $\E(G)$.

Of special interest is the left regular representation $\lambda$ of
$G$ on $\H_G = l^2(\fgeeplus)$.
This space has an orthonormal basis $\{ \xi_w : w \in \fgeeplus \}$.
We define
\[
 \lambda(p)\xi_w = \begin{cases}
                   \xi_{pw} &\qif s(p) = r(w)\\
                   0        &\quad\text{otherwise}
                   \end{cases}
\]
We will write $L_p$ for $\lambda(p)$ when $p$ is a path of positive length
and $P_x$ or $L_x$ for $\lambda(x)$ when $x \in \V(G)$ as is convenient.
Each edge $e$ is sent to a partial isometry with domain
$P_{s(e)}\H_G$ and has range contained in $P_{r(e)}\H_G$.
Moreover the ranges of $L_e$ are pairwise orthogonal for $e \in \E(G)$.
Therefore $\lambda$ is a completely contractive representation.

The norm closed algebra generated by $\lambda(\fgeeplus)$ is the \textit{tensor algebra}
of $G$ and is denoted
as $\T^{+}(G)$. The tensor algebras associated with graphs were introduced under
the name \textit{quiver algebras} by Muhly and Solel in
\cite{M,MS}. The algebras $\T^{+}(G)$ generalize Popescu's
non-commutative disc algebras \cite{Pop3} and they
are the motivating example in the theory of tensor algebras of Hilbert bimodules
\cite{M,MS}.

The free semigroupoid algebra $\fL_G$ is the \wot-closed algebra
generated by $\lambda(\fgeeplus)$ introduced by Kribs and Power \cite{KrP}.
This generalizes the free semigroup algebras $\fL_n$
which correspond to the graphs with one vertex and $n$ loop edges
\cite{Pop1, DP1}.
Each element $A$ of $\fL_G$ has a Fourier series representation
$\sum_{p \in \fgeeplus} a_p L_p$ defined as follows.\vspace{.3ex}
For each path $p \in \fgeeplus$, let $a_p = \ip{A \xi_{s(p)}, \xi_p}$.
\vspace{.3ex} Then the Cesaro means
\[
 \Sigma_k(A) = \sum_{|p|<k} \Big( 1 - \frac{|p|}k \Big) a_p L_p
\]
converge in the strong operator topology to $A$. If, in addition, $A \in \T^{+}(G)$,
then the Cesaro sums converge in norm.

A nest algebra is the algebra of all operators leaving invariant a nest, i.e.\
a complete chain of subspaces including $\{0\}$ and $\H$.
It is said to have multiplicity one if the von Neumann algebra generated by the
projections onto these subspaces is maximal abelian.
An atom of a nest is the difference of an element of the nest and its
immediate predecessor. A nest is atomic if the Hilbert space is
spanned by all of the atoms. Such nests are multiplicity free
precisely when the atoms are one-dimensional.
The finite multiplicity free nest algebras are just the upper triangular
matrices $\T_n$ acting on the $n$-dimensional Hilbert space $\bC^n$.
A type $\bN$ nest will indicate that there is an orthonormal basis
$\{h_i : i \ge 1 \}$ with subspaces $\N_k = \spn \{ h_i : i \le k \}$ for
$k \ge 0$ and $\H$.
We refer to \cite{D_nest} for more information about nest algebras.

A nest representation of an algebra $\fA$ is a representation $\rho$ of
$\fA$ on $\H$ such that $\rho(\fA)$ is \wot-dense in a nest algebra.
Our terminology differs somewhat from the original definition of
Lamoureux \cite{Lam} where he asks only that the invariant subspaces
of $\rho(\fA)$ be a nest, and not that the algebra be \wot-dense in the
nest algebra.  In the case where the \wot-closure of $\rho(\fA)$
contains a masa, these two notions are equivalent.

\section{Partially isometric dilations}

In this section we describe a dilation theory for $n$-tuples of operators which
generalizes the Frazho-Bunce-Popescu dilation Theorem
\cite{Fr1,Bun,Pop_diln} and is useful in producing representations for $\flgee$.
The results below can be found in two recent works. Muhly and Solel \cite{MS}
have obtained a dilation theory for covariant representations
of C*-correspondences which implies Theorem~\ref{dilation} below.
Jury and Kribs \cite{JK2} obtained self-contained proofs for dilations
of $n$-tuples of operators associated with directed graphs.
We follow their approach only because it avoids the technical language of
Hilbert bimodules. However we reformulate it to better reflect our interest in
representations of graphs.

Let $G$ be a countable directed graph.
A representation $\sigma$ of $\fgeeplus$ will be called \textit{partially
isometric} if $S_e :=\sigma(e)$ is a partial isometry with domain
$\sigma(s(e))\H$ for each edge $e \in \E(G)$.
We will write $P_x := \sigma(x)$ for the vertex projections.
It is easy to see that such a representation will satisfy the relations
$(\dagger)$ below, where $(1)$ is a consequence of being a
representation, $(3)$ is the partial isometry condition, and $(2)$ and
$(4)$ follow because $\rho$ is completely contractive.
\[
(\dagger)  \left\{
\begin{array}{lll}
(1)  & P_x P_y = 0 & \mbox{for all $x,y \in \V(G)$, $ x \neq y$}  \\
(2) & S_{e}^{*}S_f = 0 & \mbox{for all $ e, f \in \E(G)$, $e \neq f $}  \\
(3) & S_{e}^{*}S_e = P_{s(e)} & \mbox{for all $e \in \E(G)$}      \\
(4)  & \sum_{r(e)=x}\, S_e S_{e}^{*} \leq P_{x} & \mbox{for all $x \in \V(G)$}
\end{array}
\right.
\]

Clearly, the left regular representation is partially isometric, but there
are many others.
If a partially isometric representation $\sigma$ satisfies
\[
 \sotsum_{\substack{e \in \E(G)\\ r(e)=x}} \sigma(e) \sigma(e)^{*} = P_x
 \qforal x \in \V(G),
\]
then we say that $\sigma$ is \textit{fully coisometric}.
  From the operator algebra perspective, fully coisometric representations
generate C*-algebras which generalize the Cuntz-Krieger algebras.

At the other extreme, we say a partially isometric representation $\sigma$ is
\textit{pure} if
\[
 \sotlim_{d \to \infty} \ \Big( \sotsum_{\substack{p\in \fgeeplus\\ |p| =d} }
  \sigma(p) \sigma(p)^* \Big) = 0
\]
Pure partially isometric representations can be obtained as follows.
For $x \in \V(G)$, define the subspace
$\H_x = \spn \{ \xi_w : w \in \fgeeplus,\ w = wx \}$.
It is clear that $\H_G = \sum_{x \in\V(G)}^\oplus \H_x$
and that each $\H_x$ reduces $\lambda$; let $\lambda_x$ denote the restriction
of $\lambda$ to $\H_x$.
Observe that
\[
 \Big( \sotsum_{e \in \E(G)} L_e L_e^*\Big) \Big|_{\H_x} = P_{\H_x} - \xi_x \xi_x^*
\]
is a proper subprojection of $P_{\H_x}$, orthogonal to the unit vector $\xi_x$.
(By $xy^*$, $x , y \in \H$, we mean the rank one operator  $xy^*(z) = \ip{ z,y } x$.)
Indeed
\[
 \sotlim_{d\to\infty} \Big( \sotsum_{|p|=d} L_p L_p^*  \Big) \Big|_{\H_x} = 0
\]
because for fixed $d$, this sum is a projection onto basis vectors from the collection
\[
 \{ \xi_w : w \in \fgeeplus,\ s(w) = x,\ |w| \ge d \} ,
\]
and these sets decrease to the empty set. Therefore $\lambda$ and its subrepresentations
$\lambda_x$ are pure.

Jury and Kribs \cite{JK2} prove the converse: that a pure partially isometric
representation $\sigma$ is equivalent to the restriction of $\lambda^{(\infty)}$
to a reducing subspace.  More precisely,
$ \sigma \simeq \sum^\oplus_{x\in V(G)} \lambda_x^{(\alpha_x)}$
\vspace{.3ex}
where the multiplicities are determined by
$\alpha_x = \rank \big( P_x  ( I - \sum_e S_e S_e^* )\big)$ for $x\in\V(G)$.
They also prove an analogue of the Wold decomposition which states that every
partially isometric representation of $G$ decomposes as a direct sum of a pure
representation and a fully coisometric representation.
We have already observed that $\lambda$ and its subrepresentations
$\lambda_x$ are pure.

Let $\rho$ be a completely contractive representation of $\fgeeplus$.
A \textit{partially isometric dilation} of $\rho$ is a partially isometric
representation $\sigma$ of $\fgeeplus$ on a Hilbert space $\K$ containing $\H$
such that $\H$ is invariant for $\sigma(\fgeeplus)^*$ and
$\sigma(p)^*|_\H = \rho(p)^*$ for every path $p \in \fgeeplus$.
Say that $\sigma$ is a \textit{minimal} partially isometric dilation if in
addition $\H$ is cyclic for $\sigma(\fgeeplus)$, namely
$\K =  \bigvee_{p\in\fgeeplus} \sigma(p) \H$.
Since $\sigma(x)$ is self-adjoint, it follows that $\sigma(x)$ commutes with
$P_\H$ for all $x \in \V(G)$.

Given a row contraction $T = \begin{bmatrix}T_1 & \dots & T_n\end{bmatrix}$
or $T = \begin{bmatrix}T_1 & T_2 & \dots\end{bmatrix}$
of operators acting on a Hilbert space $\H$, consider a family of
mutually orthogonal projections $\P = \{ P_x \}_{x\in\J}$ on $\H$
which sum to the identity operator and  {\it stabilize} $T$ in the
following sense:
\[
  T_i P_x ,\, P_x T_i \in \{ T_i, 0 \} \qforal i \qandal x\in\J.
\]
Observe that these relations lead to a directed graph $G$ with
vertex set $\V(G) = \J$ and directed edges $e_i$, where $r(e_i)$ and
$s(e_i)$ are vertices with $P_{r(e_i)} T_i P_{s(e_i)} = T_i$.
Indeed these vertices are unique when $T_i \ne 0$, but are arbitrary if
$T_i=0$. We shall fix a choice of $G$ and rewrite the indices as edges of $G$,
$T = (T_e )_{e\in\E(G)}$. It is clear that this determines a unique completely
contractive representation of $\fgeeplus$ with
$\rho(e) = T_e$ and $\rho(x) = P_x$.

Let $T = \rowt$ be a row contraction stabilized by $\P = \{P_x\}_{x\in\V(G)}$
and indexed as above by a directed graph $G$.
A \textit{partially isometric $G$-dilation} of $T$ is a partially isometric
dilation $\sigma$ of the induced representation $\rho$ of $G$ determined by
$T$.

With this latter viewpoint, it was shown in \cite[Theorem 3.2]{JK2} that if $G$
has no sinks, then $T$ has a unique minimal partially isometric $G$-dilation
$\sigma$. They also show that $\sigma$ is fully coisometric if and only if $T$ is
a coisometry; that is,
\[
 I = TT^* = \sum_{e \in \E(G)} T_e T_e^* .
\]
Moreover, they show that $\sigma$ is pure if and only if
\[
 \sotlim_{d\to\infty}\ \Big( \sotsum_{\substack{p\in\fgeeplus\\|p| =d}}
  \rho(p) \rho(p)^* \Big) = 0
\]
In particular, this holds when $\|T\| < 1$.

Observe that when $\sigma$ is pure, there is a \wot-continuous extension of
$\sigma$ to $\fL_G$.  This can then be compressed to $\H$ to obtain a
\wot-continuous representation of $\fL_G$ extending $\rho$.
This provides an $\fL_G$ functional calculus for the row contraction $T$.

We require this dilation result for graphs that may have sinks.  Using a
technique of \cite{BPRS} of adding a tail to each sink, we are able to extend
the Jury--Kribs result to arbitrary countable graphs with little effort.

\begin{thm} \label{dilation}
Let $G$ be a directed graph and let $\rho$ be a completely contractive
representation on a Hilbert space $\H$. Then $\rho$ has a minimal partially
isometric dilation $\sigma$ which is unique up to unitary equivalence fixing
$\H$.
\end{thm}

\Prf The trick of \cite{BPRS} is to add an infinite tail to each sink.
Briefly, if $v_0$ is a sink, add vertices $v_i$ for $i \ge 1$ and edges $e_i$
with $s(e_i) = e_{i-1}$ and $r(e_i) = v_i$ for $i \ge 1$.
In this manner, we can embed $G$ into a larger graph $H$ with no sinks.
Extend $\rho$ to a representation $\bar{\rho}$ of $\bF^+(H)$ by
setting
$\bar{\rho}(e_i) = \bar{\rho}(v_i) = 0$ for $i \ge 1$ on every tail.
It is evident that this does not affect the condition of complete
contractivity, and $\bar{\rho}$ is pure or fully coisometric exactly when $\rho$
is. Moreover $\fgeeplus$ is a subsemigroupoid of $\bF^+(H)$ and $\fL_G$ is a
\wot-closed subalgebra of $\fL_H$.

By the Jury--Kribs dilation, there is a unique minimal partially isometric
dilation $\tau$ of $\bar{\rho}$ on a Hilbert space $\bar{\K}$ containing $\H$.
Let $\K = \bigvee \tau(\fgeeplus) \H$ and let $\sigma(p) = \tau(p)|_\K$ for all
paths $p \in \fgeeplus$.
It is easy to verify that $\sigma$ is a partially isometric dilation of
$\rho$, and minimality is easy to verify as well.
\bx

The following corollary is immediate from earlier comments.

\begin{cor}\label{dilationcor}
Let $G$ be a directed graph and let $\rho$ be a completely contractive
representation on a Hilbert space $\H$. Then the minimal partially
isometric dilation $\sigma$ of $\rho$ is fully coisometric if and only if
the row operator $T$ with entries $(\rho(e))_{e \in \E(G)}$ is a
coisometry. And $\sigma$ is pure if and only if
\[
 \sotlim_{d\to\infty}\  \Big( \sotsum_{\substack{p\in\fgeeplus\\|p| =d}}
  \rho(p) \rho(p)^* \Big) = 0 .
\]
In particular, this holds if $\|T\|<1$.
When $\rho$ is pure, it extends uniquely to a \wot-continuous completely
contractive representation of $\fL_G$.
\end{cor}

We give one routine result about the tensor algebra case in light of
Remark~\ref{R:tensor}.  This will show that the assumption of a
bounded representation of $\T^+(G)$ with a particular property implies the
existence of a \wot-continuous completely contractive representation
of $\fL_G$ with the same property.  In the converse direction, we will
always be constructing \wot-continuous completely contractive
representations of $\fL_G$.
In particular, this applies to the property of irreducibility,
being a nest representation and being faithful.
The similarity used does not affect any of these.
If we are concerned only with contractive representations,
then the restriction to the diagonal is already a $*$-representation
and that part of the proof is unnecessary.

\begin{prop}\label{P:tensor}
Suppose that $\rho$ is a bounded, norm continuous representation of
the tensor algebra $\T^+(G)$ on a Hilbert space $\H$.  Then there is a
representation similar to $\rho$ and a \wot-continuous completely
contractive representation $\phi$ of $\fL_G$ on the same space such that
the \wot-closure of the ranges of the two representations coincide.
\end{prop}

\Prf Note that the representation of $\rho$ restricted to the diagonal
algebra $\spn\{ P_x : x \in \V(G) \}$ is similar to a $*$-representation.
So we adjust $\rho$ by a similarity to achieve this.
Then there is a bound
\[
  C := \sup \{ \|\rho(e)\| : e \in \E(G) \} .
\]
Define an injective homomorphism $\tau$ of $\fL_G$ into $\T^+(G)$ by
setting $\tau(P_x) = P_x$ for $x \in \V(G)$ and
$\tau(L_{e_i}) = 2^{-i-1}C^{-1} L_{e_i}$ for any enumeration
$\{e_i : i \ge 1 \}$ of $\E(G)$.
Then Corollary~\ref{dilationcor} shows that $\rho\tau$ extends to a
\wot-continuous completely contractive representation $\phi$ of $\fL_G$.
Clearly this has the desired properties.
\bx

\section{Finite dimensional representations for free semigroupoid algebras}

First we will describe a class of representations of $\fL_G$ onto full matrix
algebras.
Let $u$ be a cycle in $G$ of length $k$ and $\lambda \in \bT$.
Write $u = e_{i_k} e_{i_{k-1}} \dots e_{i_1}$.
Let $\H_k$ be a $k$-dimensional Hilbert space with orthonormal basis
$h_1, h_2, \dots , h_k$.
(Recall that by $xy^*$, we mean the rank one operator  $xy^*(z) = \ip{ z,y } x$.)
Define a representation $\phi_{u,\lambda}$ of $\fgeeplus$ on $\H_k$ as
follows.
Set
\begin{align*}
  \phi_{u,\lambda}(x) &= \sum_{s(e_{i_j})=x} h_j h_j^*  \qand\\
  \phi_{u,\lambda}(e) &= \frac 1 2 \sum_{e_{i_j}=e}h_{j+1}h_j^*
\end{align*}
where $h_{k+1}$ denotes $\lambda h_1$.
Since the row operator with entries
\[
 \big[ \phi_{u,\lambda}(e) \big]_{e \in \E(G)}
\]
has norm $1/2$, it follows from Corollary~\ref{dilationcor} that this
extends to a completely contractive \wot-continuous representation of $\fL_G$.

\begin{defn}
A cycle $u \in \fgeeplus$ is \textit{primitive} if it is not a nontrivial
power of another cycle.
\end{defn}

Clearly every cycle factors uniquely as $u=v^p$ where $v$ is a primitive
cycle. A primitive cycle $v$ has no rotational symmetry, meaning that $v$ is
not invariant under any cyclic permutation.  Indeed we have:

\begin{lem}\label{cyclic_symmetry}
If $u$ is a cycle which can be non-trivially factored as $u= u_1 u_2 = u_2 u_1$,
then it is not primitive.
\end{lem}

\Prf The hypothesis shows that the cyclic rotation of the letters in $u$
by $|u_1|$ steps leaves $u$ invariant.
The set of cyclic permutations leaving the cycle $u$ invariant is a
subgroup of the cyclic group $C_{|u|}$, and hence is $C_d$ for some divisor
$d$ of $|u|$.  Since this subgroup is non-trivial, we obtain $u = v^d$ where
$v$ is the cycle formed by the first $|u|/d$ edges of $u$.  In particular, $u$
is not primitive.
\bx

If $u=v^p$, then the representation $\phi_{u,\lambda}$ will not be irreducible.
To see this, define a unitary $U$ by $U h_j = h_{j + k/p}$ where we identify
$h_{k+j}$ with $\lambda h_j$.
Observe that $U$ commutes with the range of $\phi_{u,\lambda}$.
We now establish the converse.

\begin{lem}  \label{irredlem}
When $u$ is a primitive cycle and $\lambda \in \bT$, the representation
$\phi_{u,\lambda}$ is irreducible.
\end{lem}

\Prf Let $u_j = e_{i_{j-1}} \dots e_{i_2} e_{i_1} e_{i_k} e_{i_{k-1}} \dots
e_{i_j}$ be the cycle that follows $u$ around beginning with the $j$th edge.
Primitivity ensures that the cycles $u_j$ are distinct.
Observe that $\phi_{u,\lambda}(u) = 2^{-k} \lambda h_1 h_1^*$.
Indeed, if you apply $\phi_{u,\lambda}(u)$ to each basis vector $h_j$, one
sees that $h_1$ returns to a multiple of itself.  But starting with any other
$h_j$, the result must be $0$ unless at each step the edge $e_{i_s}$
coincides with the edge $e_{i_{s+j}}$, which cannot occur for all $k$ edges
by primitivity.
Similarly $\phi_{u,\lambda}(u_j) = 2^{-k} \lambda h_j h_j^*$.

Hence we see that
$\phi_{u,\lambda}(e_{i_j}u_j) = 2^{-k-1} \lambda h_{j+1} h_j^*$.
So the range of $\phi_{u,\lambda}$ contains a set of matrix units for
$\B(\H_k)$, whence it is irreducible.
\bx

Recall that $G$ is transitive if there is a path from each vertex in
the graph to every other.
It is not hard to see that each connected component of $G$ is transitive
if and only if every  edge of $G$ lies on a cycle.

\begin{thm} \label{irreducible}
If $G$ is a countable directed graph, then the finite dimensional irreducible
representations of $\flgee$
separate points if and only if $G$ is transitive in each component.
\end{thm}

\Prf  If the finite dimensional irreducible
representations of $\flgee$ separate the points, then $\flgee$ is semisimple
and therefore $G$ is transitive in each component \cite[Theorem 5.1]{KrP}.

Conversely, assume that $G$ is transitive in each component.
For any path $w\in \fgeeplus$, the transitivity in each component of $G$
implies that there exists a path $v \in \fgeeplus$ so that $vw$ is a cycle.
Hence there is a primitive cycle $u = e_{i_k} e_{i_{k-1}} \dots e_{i_1}$
so that $vw = u^p$. Consequently there exist $0 \le k'<k$ and $n\ge0$ so that
$w = e_{i_{k'}} \dots e_{i_1} u^n$, where we mean $u^n$ if $k'=0$.
Let us write $p_r = e_{i_{k'}} \dots e_{i_1} u^r$ for $r \ge 0$.

We will make use of the irreducible representations $\phi_{u,\lambda}$ for
$\lambda \in \bT$.
Let us evaluate $g(p,\lambda) := \ip{ \phi_{u,\lambda}(p) h_1, h_{k'} }$ for
paths $p \in \fgeeplus$.
Arguing as in Lemma~\ref{irredlem}, one sees that $g(p,\lambda) = 0$ unless
$p = p_r$ for some $r\ge0$; and $g(p_r) = \lambda^r 2^{-k'-rk}$.
This power of $\lambda$ will allow us to distinguish the $p_r$ by integration.

For each $A \in \flgee$, let
$g(A,\lambda) = \ip{ \phi_{u,\lambda}(A) h_1, h_{k'} }$.
If $A$ has the Fourier series $\sum_{p \in \fgeeplus} a_p L_p$,
then,
\begin{equation} \label{Cseries}
  g(A,\lambda) = \sum_{r \ge0} 2^{-k'-rk} a_{p_r} \lambda^r .
\end{equation}
The above equality requires a short argument. The calculation of the
previous paragraph shows that
$g(A,\lambda)$ equals the limit of the Cesaro partial sums of the series in display.
However, the coefficients $2^{-k'-rk}a_p$ are dominated by a geometric series. This
is easily seen to imply that the Cesaro partial sums as well as the series in
(\ref{Cseries}) converge
absolutely and uniformly over the circle $\bT$ to the same limit, i.e., $g(A,\lambda)$.
We may now compute
\[
 a_w = \frac{2^{k'+rn}}{2\pi i} \int_\bT \bar{\lambda}^{n+1} g(A,\lambda)
 \,d\lambda.
\]
For any non-zero element $A \in \fL_G$ , there exists a least one non-zero
Fourier coefficient $a_w$ for some path $w$.
The argument above shows that $g(A,\lambda)$ is a non-zero function, and hence
$\phi_{u,\lambda}(A) \ne 0$ for some $\lambda\in\bT$.
Therefore these representations separate points.
\bx

Recall that the strong radical of a Banach algebra is the intersection of all
maximal ideals.  If this intersection is $\{0\}$, the algebra is called
strongly semisimple.
The following corollary is new even for the free semigroup algebras $\fL_k$,
$k\geq2$.

\begin{cor} \label{strongrad}
The following are equivalent:
\begin{enumerate}
\item $\flgee$ is semisimple
\item $\flgee$ is strongly semisimple
\item $G$ is transitive in each component.
\end{enumerate}
If $\V(G)$ is finite, then the Jacobson radical of $\flgee$ is equal to the
strong radical.
\end{cor}

\Prf The kernels of finite dimensional irreducible representations are maximal
ideals.  Thus (3) implies (2) by Theorem~\ref{irreducible}.
Clearly (2) implies (1).
By Kribs and Power \cite[Theorem 5.1]{KrP}, (1) implies (3).

If $\V(G)$ is finite, Kribs and Power also show that the Jacobson radical is
generated by $\lambda(e)$ for those edges $e\in\E(G)$ which are not contained
in any cycle.
Indeed, a little thought shows that this ideal is the \wot-closed span of
$\lambda(w)$ of all paths $w \in \fgeeplus$ which are not contained in any
cycle (in the sense that there is no path $v$ such that $vw$ is a cycle).
This means that $A \in \fL_G$ lies in the radical if and only if the Fourier
coefficients $a_w$ vanish whenever $w$ is contained in a cycle.
On the other hand, suppose that $A \in \fL_G$ and $a_w \ne 0$ for some
$w\in\fgeeplus$ which is contained in a cycle $vw$.
Then factor $vw=u^k$ as a power of a primitive cycle $u$.
The proof of Theorem~\ref{irreducible} shows that $A$ is not in the
maximal ideal $\ker \phi_{u,\lambda}$ for some values of $\lambda \in \bT$;
whence $A$ is not in the strong radical.  So the two radicals coincide.
\bx

The previous theorem excludes a variety of interesting graphs. In order to
cover more cases, we replace irreducible representations with nest
representations, where we seek representations that map $\fL_G$ onto a
\wot-dense subalgebra of a nest algebra.  When restricting attention to
finite dimensional representations, the nest algebras are the set of block
upper triangular operators (with respect to some decomposition of a finite
dimensional space into orthogonal subspaces).  Moreover, density implies
surjectivity in this case; and the \wot-density of the span of $\fgeeplus$
means that its image spans the nest algebra.
Finite dimensional irreducible representations are a special case because $\M_n$
is a nest algebra.

\begin{lem}  \label{decomp}
Let $G$ be a countable directed graph.
Any path $w \in \fgeeplus$ decomposes uniquely as a product
\[
  w = w_l v_l w_{l-1} v_{l-1} \dots w_1 v_1 w_0
\]
where all the edges involved in a path $w_i$ lie in a single transitive
component of $G$ $($or the term is a vertex$)$ and each $v_i$ is an edge
with source and range in distinct transitive components.
\end{lem}

\Prf If there is
an edge from one transitive component to another, there can be no path back.
The directed graph obtained by mapping each transitive component to a single
vertex has no (directed) cycles.
It follows that $s(v_i)$ lie in distinct transitive components, as do $r(v_i)$.
If there is no path $w_i$, set it to be the vertex $w_i = r(v_i) = s(v_{i+1})$.
\bx

\begin{thm} \label{fdnest_separate}
If $G$ is a countable directed graph, then the finite dimensional nest
representations separate the points in $\flgee$.
\end{thm}

\Prf
Once again, we will construct a family of nest representations which allow us
to recover the Fourier coefficients in $\flgee$.

Let $w$ be a path in $\fgeeplus$.
Factor $w$ using Lemma~\ref{decomp} above as
$w = w_l v_l w_{l-1} v_{l-1} \dots w_1 v_1 w_0$.
Since each $w_i$ is contained in a single transitive component, arbitrarily
complete it to a cycle $w'_i w_i$, which we then factor as a power of a primitive
cycle, say $w'_i w_i = u_i^{k_i}$ and $w_i = p_i u_i^{n_i}$ for some proper
subpath $p_i$ of $u_i$.   If $w_i$ is a vertex, set $u_i = w_i$.
Let $\phi_{u_i,\lambda_i}$ be the irreducible representations defined at the
beginning of this section, each acting on a finite dimensional Hilbert space $\H_i$.
If $u_i$ is a vertex, take $\H_i$ to be one-dimensional and set
$\phi_{u_i,\lambda_i}(u_i)=1$, and $\phi_{u_i,\lambda_i}(p) = 0$ for all other
paths $p \in \fgeeplus$.

Set $\H = \sum^\oplus_i \H_i$. For each path $w\in\fgeeplus$ and
$\lambda = (\lambda_i)_{i=0}^l \in \bT^{l+1}$,
define a representation $\rho_{w,\lambda}$ as follows.
First we set
\[
  P_{\H_i}\rho_{w,\lambda}|_{\H_i} = \phi_{u_i,\lambda_i} \qfor 0 \le i \le l
\]
Observe that $\phi_{u_i,\lambda_i}$ is zero except on paths which lie entirely in
the transitive component of $w_i$.  In particular, this annihilates all edges
except those with both source and range in this transitive component.
Now $\phi_{u_i,\lambda_i}(w_i)$ has the form $2^{-|w_i|} \lambda_i^{n_i}
k_i h_i^*$ for certain vectors $k_i, h_i \in \H_i$.
Set $\rho_{w,\lambda}(v_i) = \frac12 h_i k_{i-1}^*$.
All other edges and vertices are sent to $0$.
The image of a path is determined by the product of the edges.
As before, we verify that the row operator with entries
$\big[ \rho_{w,\lambda}(e) \big]_{e\in\E(G)}$ has norm $1/2$, and thus by
Corollary~\ref{dilationcor} this representation extends from $\fgeeplus$ to
$\fL_G$ to be a completely contractive \wot-continuous representation.

To see that this is a nest representation, observe that the proof of
Theorem~\ref{irreducible} shows that the blocks $\B(\H_i)$ lie in the image.
The only other non-zero edges are $\rho_{w,\lambda}(v_i) = \frac12 h_i k_{i-1}^*$.
These provide the connection from $\H_{i-1}$ to $\H_i$.  It is evident that
each element is block lower triangular, and that these operators generate the
block lower triangular nest algebra as an algebra.

Now consider the function $g(A,\lambda) = \ip{ \rho_{w,\lambda}(A) h_0, k_l }$
for $A \in \fL_G$ and $\lambda \in \bT^{l+1}$.
The idea now is much as before.
We observe that if $p$ is a path, then $g(p,\lambda) = 0$ unless
\[
  p_m := p_l u_l^{m_l} v_l p_{l-1} u_{l-1}^{m_{l-1}} \dots v_1 p_0 u_0^{m_0}
\]
for $m = (m_0,\dots, m_l) \in \bZ^{l+1}$;
and $g(p_m,\lambda) = 2^{-|p_m|} \prod_{i=0}^l \lambda_i^{m_i}$.
If $A$ has Fourier series $A \sim \sum_{p \in \fgeeplus} a_p L_p$, then
\[
 g(A,\lambda) =
 \sum_{m \in \bZ^{l+1}} a_{p_m} 2^{-|p_m|} \prod_{i=0}^l \lambda_i^{m_i} .
\]
This series converges absolutely and uniformly on a neighbourhood of $\bT^{l+1}$ to
a holomorphic function because the Fourier coefficients are bounded by $\|A\|$.
Thus we recover the coefficients by integration:
\[
 a_w = 2^{|w|} \frac1{(2\pi i)^{l+1}} \int_{\bT^{l+1}} g(A,\lambda)
               \prod_{i=0}^l \bar{\lambda}_i^{n_i+1 } \, d\lambda_i .
\]
Therefore there exists some word $w$ and $\lambda$ so that
$\rho_{w,\lambda}(A) \ne 0$,

We conclude that the collection of all nest representations separates points.
\bx

Since the functions $g(A,\lambda)$ are continuous, it suffices to use a
countable set of $\lambda$'s for each (of countably many) words
$w \in \fgeeplus$ in order to separate points.

\begin{cor} \label{cor:Fourier}
Let $G$ be a countable directed graph. Then is a countable family of nest
representations that separates the points in $\flgee$.
\end{cor}

A similar conclusion about the countability of a separating family of
irreducible representations holds when $G$ is transitive in each component
and in the case to follow of upper triangular nest representations.

An ideal is meet irreducible if it is not the intersection of two strictly larger
ideals.  The kernel of a nest representation is always meet irreducible.
In particular in the finite dimensional case, this holds because a nest algebra on a
finite dimensional space has a unique minimal ideal.
We can now state the appropriate generalization of Corollary \ref{strongrad} to
arbitrary free semigroupoid algebras.

\begin{cor}
The meet irreducible ideals of a free semigroupoid algebra $\flgee$ have zero
intersection.
\end{cor}

Finally in this section, we consider representations onto the upper triangular
matrices $\T_n$ with respect to a fixed basis $h_1,\dots, h_n$ of $\H_n$.
We call such representation an \textit{upper triangular nest representation}.
This constraint puts special requirements on the graph. Note that the result below is
essential in the proof of Theorem \ref{Nrepn}.

\begin{thm} \label{upper}
If $G$ is a countable directed graph, then the upper triangular nest
representations of $\flgee$ separate points if and only if every vertex $x
\in \V(G)$ which supports a cycle also supports at least one loop edge.
\end{thm}

\Prf Suppose that there exists a cycle $w$ with $r(w)=s(w)=x$ but $x$ supports
no loop edges. We claim that any upper triangular nest representation $\psi$ of
$\fgeeplus$ satisfies $\psi(w)=0$.

Since $\psi(x)$ is a projection in $\T_n$, it is diagonal.
We will show that it always has rank at most one.
Notice that $x \fgeeplus x $ consists of all cycles $u\in \fgeeplus$ with
$r(u)=s(u)=x$.
Since $x$ supports no loop edges, this cycle must pass through another vertex
$y$.  Therefore we can factor $u$ as $u = x u_1 y u_2 x$.
Hence
\[
 \psi(u) = \psi(x) \psi(u_1) \psi(y) \psi(u_2) \psi(x) .
\]
Since $\psi(x)$ and $\psi(y)$ are orthogonal diagonal projections, this means
that $\psi(u)$ is strictly upper triangular.
Let $\Delta_n$ denote the expectation of $\M_n$ onto the diagonal $\D_n$,
and recall that the restriction of $\Delta_n$ to $\T_n$ is a homomorphism.
So
\begin{align*}
  \psi(x) \D_n &= \psi(x) \Delta_n(\spn \psi(\fgeeplus)) \psi(x) \\ &=
  \spn \Delta_n( \psi(x \fgeeplus x) ) = \bC \psi(x)
\end{align*}
because $\Delta_n(\psi(u)) = 0$ all non-trivial cycles $u$ beginning at $x$.
Therefore $\psi(x) \D_n$ is one dimensional; whence so is $\psi(x)$.
Consequently $\psi(x) \T_n \psi(x) = \bC \psi(x)$ supports no strictly upper
triangular operators.  So $\psi(u) = 0$ for all cycles $u$ beginning at $x$.
Thus these representations cannot separate points.

Conversely, suppose that every vertex $x\in\V(G)$ which supports a cycle also
supports a loop edge.  We proceed as in the proof of Theorem~\ref{irreducible}
to construct a family of upper triangular nest representations of $\fL_G$ from
which we can recover the Fourier coefficients of each element $A \in \fL_G$.

Let $\L$ denote the set of vertices $x\in\V(G)$ which supports a loop edge.
For each vertex $x \in \L$, select one loop edge at $x$ and label it $f_x$.
Call the remaining vertices $\R = \V(G) \bsl \L$ the \textit{rank one} vertices.

Let $w = e_{i_{k-1}} e_{i_{k-2}} \dots e_{i_2} e_{i_1}$ be a path in
$\fgeeplus$ which does not involve any of the designated loop edges.
Denote the vertices involved by $x_j = s(e_{i_j})$ for $1 \le j < k$ and
$x_k = r(e_{i_{k-1}})$. Let $J = \{ j : x_j \in \L \}$.
For each $\lambda = (\lambda_j)_{j \in J} \in \bT^J$ with \textit{distinct} coefficients,
we will define a representation $\psi_{w,\lambda}$ of $\fgeeplus$ on $\H_k$.
Set
{\allowdisplaybreaks
\begin{alignat*}{2}
 \psi_{w,\lambda}(x)   &= \sum_{x_j=x} h_j h_j^*
       &\FOR& x \in \V(G) \\
 \psi_{w,\lambda}(f_x) &= \frac 1 2 \sum_{x_j=x} \lambda_j h_j h_j^*
       &\FOR& x \in \L\\
 \psi_{w,\lambda}(e)   &= \frac 1 2 \sum_{e_{i_j}=e}h_{j+1}h_j^* \quad
       &\FOR& e \in \E(G) \bsl \{f_x : x \in \L\}
\end{alignat*}
}  

The row operator $E$ with entries $(\psi_{w,\lambda}(e))_{e \in \E(G)}$
satisfies
\[
 EE^* = \frac 1 4 \big( (I - h_1 h_1^*) + \sum_{j \in J} h_j h_j^* \big)
 \le \frac 1 2 I .
\]
This is a strict contraction, and thus by Corollary~\ref{dilationcor}, the
representation $\psi_{w,\lambda}$ extends uniquely to a completely contractive
\wot-continuous representation of $\fL_G$.

It is clear that the range of $\psi_{w,\lambda}$ is contained in
the algebra $\T_n^*$ of lower triangular matrices.
We will show that it is surjective.
If $j \not\in J$, then $x_j$ is a rank one vertex.
Since $x$ does not lie on any cycle, it can occur only once on the path $w$.
Therefore $\psi_{w,\lambda}(x_j) = h_j h_j^*$.
Now consider a vertex $x$ such that $x=x_j$ for some $j \in J$ and let
$J_x = \{ j \in J : x_j = x \}$.
Then
\[
 \psi_{w,\lambda}(2f_x) = \sum_{j \in J_x} \lambda_j h_j h_j^*
\]
is a diagonal operator with distinct eigenvalues $\lambda_j$.
Therefore the algebra it generates contains each $h_j h_j^*$,
Thus the image of $\psi_{w,\lambda}$ contains the diagonal.
Moreover $\psi_{w,\lambda}(e_{i_j}) h_j h_j^* = h_{j+1} h_j^*$.
Therefore the range of $\psi_{w,\lambda}$, which is an algebra, is all of
$\T_n^*$.
So by reversing the order of the basis, we obtain an upper triangular
nest representation.

Now let $v\in \fgeeplus$ be given.
We distinguish between edges of the form $f_x$ and the rest.
Deleting all occurrences of the $f_x$'s leaves a path $w$ as above.
Then putting these loops back yields a factorization of the form
\[
 v = f_{x_k}^{m_k} e_{i_{k-1}} f_{x_{k-1}}^{m_{k-1}} e_{i_{k-2}} \dots
                   e_{i_2} f_{x_2}^{m_2} e_{i_1} f_{x_1}^{m_1}
   =: p_{w,m}
\]
where for notational convenience we allow only $m_j=0$ when $j \not\in J$
and ignore the fact that $f_x$ does not exist in this case.

We now proceed as in the proof of Theorem \ref{irreducible}.
For each path $w$ above and for each $A \in \flgee$, define a function almost
everywhere on $\bT^J$ (i.e. when the coefficients are distinct) by
\[
 f_w(A,\lambda) = \ip{ \psi_{w,\lambda}(A) h_1, h_k } .
\]
Observe that as in the earlier proof that for a path $p$, one will have
$f_w(p,\lambda) = 0$ unless $p$ has the form $p_{w,m}$ for some
$m=(m_1,\dots,m_k)$; and
\[
 f_w(p_{w,m}) = 2^{-|p_{w,m}|} \prod_{j \in J} \lambda_j^{m_j}
\]
where $|p_{w,m}| = (k-1) + |m|$.  (Here $|m| = \sum m_j$.)
Therefore if $A$ has Fourier expansion $A \sim \sum_{p \in \fgeeplus} a_p L_p$,
then
\[
 f_w(A,\lambda) = \sum_{m \in \bN^J} 2^{1-k-|m|}
                   a_{p_{w,m}} \prod_{j \in J} \lambda_j^{m_j} .
\]
As in the previous proof, observe that this series converges absolutely and
uniformly on a neighbourhood of $\bT^J$ to a holomorphic function because the
Fourier coefficients are bounded by $\|A\|$. Thus we can recover
\[
 a_v = a_{p_{w,m}}
     = 2^{k-1 + |m|} \frac1{(2\pi i)^{|J|}} \int_{\bT^J} f_w(A,\lambda)
                     \prod_{j \in J} \bar{\lambda}_j^{m_j+1 } \, d\lambda_j
\]
In particular, if $a_v \ne 0$, then there is a set of positive measure on which
$f_w(A,\lambda)$ is non-zero.  So there are representations with
$\psi_{w,\lambda}(A) \ne 0$.  Thus they separate points.
\bx

\section{Faithful nest representations}

In this section, we consider the question of when a free semigroupoid algebra
$\fL_G$ has a faithful nest representation.  The first case to consider is the
nest algebra $\B(\H)$---in other words, a faithful irreducible representation.

It is useful to first describe some faithful irreducible representations of the
free semigroup algebras $\fL_n$.  The argument of Read \cite{Read} provides a
much stronger statement, namely that there is a $*$-representation of $\O_2$
which carries the non-commutative disk algebra onto an algebra which is weak-$*$
dense in $\B(\H)$.  See \cite{D_read} for a simpler proof which works for all
$n \ge 2$.  However these examples are not so easy to work into our context;
and indeed, our purpose may be achieved much more easily anyhow.

This first lemma is likely well-known, but the best we could find in the
literature doesn't quite do it, although it comes close.
Indeed in \cite[Corollary~1 of Theorem~3.6]{FW}, they construct a unitary
group $(U_t)_{t \in \bR}$ for a given non-closed  operator range $A\H$ so
that $\ran(U_t A U_t^*)$ are pairwise disjoint.

\begin{lem}\label{op_range}
There is an injective compact operator
\[
  K = \begin{bmatrix}K_1 & K_2 & \dots\end{bmatrix}
\]
mapping $\H^{(\infty)}$ into $\H$ so that the range of each $K_i$ is
dense in $\H$.
\end{lem}

\Prf We follow the proof of \cite[Theorem~3.6]{FW}.
They work with $\H = L^2(0,2\pi)$ and begin with $A = \diag (e^{-n^2})$ with
respect to the basis $e_n = e^{inx}$ for $n \in \bZ$.
They observe that the rapid decrease of the Fourier coefficients of functions
in the range means that each agrees with an entire function restricted
to $[0,2\pi]$ (almost everywhere, of course).
Then they define a unitary $U$ as multiplication by the function
$\upchi_{(0,\pi)} - \upchi_{(\pi, 2\pi)}$, and observe that conjugating $A$
by $U$ yields an operator range with trivial intersection.

For $k\ge1$, let $U_k$ be the multiplication operator on $L^2(0,2\pi)$
by
\[
  \upchi_k := \upchi_{(0,2\pi/k)} - \upchi_{(2\pi/k, 2\pi)} .
\]
Let $K_k = 2^{-k} U_k A U_k^*$ for $k \ge 1$ and
$K = \begin{bmatrix}K_1 & K_2 & \dots\end{bmatrix}$.

A function in the range of $K$ is an absolutely convergent sum of functions
which will then be the restrictions of entire functions on each interval
$[2\pi/(k+1), 2\pi/k]$ for each $k \ge 1$.
The change at the point $2\pi/k$ depends only on the image of $K_k$ for $k
\ge2$.
Suppose that a vector $x = \sumoplus_{k\ge1}x_k$ satisfies $Kx = 0$.
As noted above, the image $K_k x_k = \upchi_k f_k$, where $f_k$ is the
restriction of an entire function.
So if $y_k$ is the vector $x$ with the $k$th coordinate replaced by $0$.
Then $g_k = Ky_k$ agrees with an entire function on $[2\pi/(k+1), 2\pi/(k-1)]$.
So $\upchi_k f_k$ must agree with $g_k$ on this interval.
This forces $f_k$ and all of its derivatives to vanish at $2\pi/k$, whence
$f_k=0$ and thus $x_k=0$.
We conclude that $x_k = 0$ for $k\ge2$; and since $K_1$ is injective, so also
is $x_1=0$.  Therefore $K$ is injective.
\bx

\begin{lem}\label{compacts}
Suppose that $K_1$ and $K_2$ are injective compact operators with dense range.
Then there are invertible operators $S_i$ so that the algebra generated by
$K_1S_1$ and $K_2S_2$ is the whole algebra $\fK$ of compact operators.
\end{lem}

\Prf Observe that the injectivity and the range of $K_i$ is unchanged if
it is replaced by $K_i S_i$ when $S_i$ is invertible.
There is a unitary $U$ so that $K_1 U$ is positive, and thus is
diagonalizable with respect to some basis $\{ e_i : i \ge 1 \}$.
Then choosing $S_1 = UD$ for some invertible diagonal operator $D$,
we may arrange that $K_1 S_1$ has distinct eigenvalues.

We wish to modify $K_2$ to $K_2 S_2$ so that the matrix entries $\ip{ K_2 S_2
e_1,e_i}$ and $\ip{K_2 S_2 e_i, e_1}$ are non-zero for all $i \ge 2$.
Since the range of $K_2$ is dense, it is routine to verify that the range
contains a vector $y = K_2 x$ so that $\alpha_i = \ip{y, e_i} \ne 0$ for
all $i \ge 1$; and in addition, we may easily arrange that $x$ is not a
multiple of $K_2^* e_1$.
Let $z = P_{\bC x}^\perp K_2^* e_1$.
Then pick a vector $w$ so that
$\ip{w,e_1}=0$ and $\ip{w,e_i} \ne 0$ for all $i \ge 2$.
Finally let $V$ be any isometry mapping $\spn\{e_1,w\}^\perp$ onto
$\spn\{x,z\}^\perp$.
We now define an invertible operator
\[
 S_2 = xe_1^* + \|z\|^{-2}zw^* + V .
\]
It now follows that for all $i \ge 1$,
\[
 \ip{K_2 S_2 e_1, e_i} = \ip{K_2 x, e_i} = \ip{y,e_i} = \alpha_i \ne 0
\]
and
\begin{align*}
 \ip{ K_2 S_2 e_i, e_1 } &= \ip{ e_i, S_2^* K_2^* e_1 }\\
 &= \ip{ e_i, S_2^* (\bar{\alpha}_1 \|x\|^{-2} x + z) } \\
 &= \ip{ e_i, \bar{\alpha}_1 e_1 + w} \ne 0 .
\end{align*}

We may replace $K_2$ by $K_2S_2$ without changing its range.
However we now have that $\Alg\{K_1, K_2\}$ is the whole algebra of compact
operators.  This is because $K_1$ generates the diagonal algebra with respect
to the basis $\{e_i : i \ge 1 \}$, and the non-zero matrix entries of $K_2$
ensure that the matrix units $e_i e_1^*$ and $e_1 e_i^*$ lie in the algebra
for all $i \ge 1$.  This guarantees that the whole set of compact operators is
generated as an algebra by $K_1$ and $K_2$.
\bx

Recall that the free semigroup algebras $\fL_n$, $2 \le n \le \infty$, are
the free semigroupoid algebras for the graph $P_n$ with one vertex and $n$
loop edges. We write $\bF^+_n$ for the free semigroup on $n$ letters, $\fL_n$
for the algebra generated by the left regular representation, and $\H_n$ for
the Hilbert space rather than $\bF^+(P_n)$, $\fL(P_n)$ and $\H_{P_n}$.

\begin{thm}\label{Ln_irred}
The free semigroup algebras $\fL_n$, for $2 \le n \le \infty$, have faithful
irreducible representations.
\end{thm}

\Prf We will use notation as if $n$ is finite, but the argument works
equally well if $n = \infty$.
Let $\H$ be a separable Hilbert space.
Use Lemma~\ref{op_range} to obtain an injective compact operator
$K = \begin{bmatrix}K_1 &\! \dots \! & K_n\end{bmatrix}$ mapping
$\H^{(n)}$ into $\H$, so that the range of each $K_i$ is dense.
Use Lemma~\ref{compacts} to replace $K_i$ by $K_i S_i$ for $i=1,2$ so that the
algebra that they generate is all of the compact operators.
Scale each $K_i$ by a constant in $(0,1]$ so that $\sum_{i \ge 1} \|K_i\| = r
< 1$.
This does not affect the properties that we have established.

In particular, $\|K\| \le r < 1$.
By \cite{Pop2} or by Corollary~\ref{dilationcor}, there is a unique completely
contractive \wot-continuous representation $\rho$ of $\fL_n$ into $\B(\H)$
such that $\rho(L_i) = K_i$ for all $i \le n$.
Since $\Alg \{ K_1, K_2 \} = \fK(\H)$, it follows that this representation is
irreducible.

Since $K$ is compact and has dense range, $\ran K$ is not closed.
Pick a unit vector $x_\mt$ which is not in this range.
Define vectors $x_w = \rho(L_w) x_\mt$ for all $w \in \bF^+_n$.
We will show that $V := \sum_{w \in \bF^+_n} x_w \xi_w^*$ is a bounded
operator from $\H_n$ into $\H$.
Indeed setting $\alpha = (\alpha_i)_{i=1}^n$ where $\alpha_i = \| K_i \|$,
\begin{align*}
  \|V\| &\le \sum_{w \in \bF^+_n} \| x_w \|
         \le \sum_{w \in \bF^+_n} \| \rho(L_w) \|
         \le \sum_{w \in \bF^+_n} w (\alpha) \\
        &\le \sum_{k \ge 0} \sum_{|w|=k} w (\alpha)
         =   \sum_{k \ge 0} \big( \sum_{i=1}^n \alpha_i \big)^k
         \le \sum_{k \ge 0}  r^k = \frac 1 {1-r} .
\end{align*}

Next we claim that $V$ is injective.
Suppose, to the contrary, that $V \zeta = 0$ where $\zeta =
\sum_{w \in \bF^+_n} a_w \xi_w$ in non-zero.
Let $w_0$ be a word of minimal length $k_0$ with $a_{w_0} \ne 0$.
We may write this vector uniquely as a sum
\[
  \zeta = \sum_{|w|=k_0} \zeta_w
        = \zeta_{w_0} + \sum_{|w|=k_0,\ w \ne w_0} \zeta_w
        =: \zeta_{w_0} + \eta
\]
where
\[
 \zeta_w = \sum_{v \in \bF_n^+} a_{wv} \xi_{wv} \qand
 \eta = \sum_{|w|=k_0,\ w \ne w_0} \zeta_w .
\]
Note that $V\zeta_w$ lies in the range of $w(K_1,K_2,\dots)$.
Observe that the row operator $J_{k_0}$ from $\H^{n^{k_0}}$ into $\H$ with
entries $\big( w(K_1,K_2,\dots) \big)_{|w|=k}$ is injective.
In particular, the range of $w_0(K_1,K_2,\dots)$ does not intersect the range
of the row operator with entries
$\big( w(K_1,K_2,\dots) \big)_{|w|=k,\ w\ne w_0}$.
In particular, $0 = V\zeta = V \zeta_{w_{0}}  + V\eta$ implies that
$V \zeta_{w_{0}}  = V\eta = 0$.

Now $\zeta_{w_{0}} = a_{w_0} \xi_{w_0} + \sum_{i \ge 1} \zeta_{w_0 i}$.
Arguing as before, $V \sum_{i \ge 1} \zeta_{w_0 i}$ lies in the range of
$w_0(K_1,K_2,\dots) \ran K$ while
\[
 a_{w_0} x_{w_0} = V a_{w_0} \xi_{w_0} = w_0(K_1,K_2,\dots) a_{w_0} x_\mt
\]
does not.
The injectivity of $w_0(K_1,K_2,\dots)$ shows that $a_{w_0} = 0$, contrary to
our assumption.
Therefore $V$ is injective.

Clearly $V$ is an intertwining map for $\rho$, meaning that $\rho(A) V = V A$
for all $A \in \fL_n$, as this holds for $\{L_w : w \in \bF_n^+ \}$ and extends
by \wot-continuity to the whole algebra.
We can now show that $\rho$ is faithful.
If it is not, then the kernel is a non-zero \wot-closed ideal of $\fL_n$.
By \cite[Theorem~2.1]{DP2}, such an ideal must contain an isometry $L$.
But then $0 = \rho(L) V = V L$ contradicts the injectivity of $V$.
So $\rho$ is faithful.
\bx

\begin{rem}\label{separating}
There is an interesting and unexpected consequence of this proof, which is
that the irreducible representations constructed have a separating vector,
namely $x_\mt$.  Indeed, if $A \in \fL_n$ is non-zero, then
\[
 \rho(A) x_\mt = \rho(A) V \xi_\mt = V A \xi_\mt .
\]
As $A \xi_\mt \ne 0$ and $V$ is injective, the result follows.
Of course $x_\mt$ is also cyclic, as is every non-zero vector.
For future reference, we note that the only condition used in the choice of
$x_\mt$ was that it was not in the range of $K$.
\end{rem}

We now extend the argument to the general case.

\begin{defn}\label{strong_transitive}
A graph $C$ is a \textit{cycle graph} if it has a finite set
$\V(C) = \{ x_1,\dots,x_n \}$ of vertices and edges
$\E(C) = \{ e_1,\dots,e_n \}$ with
$s(e_i) = x_i$ and $r(e_i) = x_{i+1 \pmod n}$.
The \textit{trivial graph} consists of one vertex and no edges.
A non-trivial transitive graph which is not a cycle is called
\textit{strongly transitive.}
\end{defn}

Note that a single vertex with one loop edge is transitive but not strongly
transitive.

\begin{thm}\label{faithful_irred}
Let $G$ be a non-trivial countable directed graph.  Then $\fL_G$ has a faithful
irreducible representation if and only if it is strongly transitive.
\end{thm}

\Prf Transitivity is an obvious necessary condition: for if $x$ and $y$ are
vertices and $\phi$ is a faithful representation, the projections $\phi(x)$ and
$\phi(y)$ are non-zero and $\phi(y) \phi(\fL_G) \phi(x)$ will be non-zero if and
only if there are paths from $x$ to $y$.
Also $G$ cannot be a cycle for then $x \fgeeplus x$ is abelian and
$P_x \fL_G P_x$ is isomorphic to $H^\infty$.
Clearly any representation will restrict to an abelian algebra on
$\rho(P_x)\H$.  If this is dense, then $\rank \rho(P_x) = 1$.
But then $\rho$ is not injective unless there are no cycles at all, which
would mean that $G$ is the trivial graph.

We will establish the converse.
Let $G$ be a strongly transitive graph.
Observe first that if each vertex in a transitive graph is the source
of exactly one edge, then $G$ is a cycle.
So there is a least one vertex $x_0$ and two distinct edges $e_1$ and $e_2$
with $s(e_1) = s(e_2) = x_0$.

Let $\E(G)$ be enumerated as $\{ e_1, e_2, \dots \}$.
Using Lemma~\ref{op_range}, we select an injective compact operator
\[
 K = \begin{bmatrix} K_1 & K_2 & \dots\end{bmatrix}
\]
so that each $K_i$ has dense range.
This property is unaffected if each $K_i$ is multiplied by non-zero
scalars bounded by 1, so we may also assume that
$\sum_{i\ge1} \|K_i\| = r < 1$.
Therefore $\|K\| \le r$.

We will define a representation $\rho$ on a Hilbert space
\[
 \K = \H^{|\V(G)|} = \sum_{x \in \V(G)} \!\!\oplus\ \H_x
\]
so that $Q_x := \rho(x) = P_{\H_x}$ for each vertex $x \in \V(G)$.
Let $J_x$ be a fixed unitary mapping $\H$ onto $\H_x$.
By Corollary~\ref{dilationcor}, there is a \wot-continuous representation
$\rho$ of $\fL_G$ on such that $\rho(e_i) = J_{r(e_i)} K_i J_{s(e_i)}^*$.

Pick two minimal cycles with source $x_0$ of the form $c_1 = p_1 e_1$ and
$c_2 = p_2 e_2$.
Observe that the edges $e_1$ and $e_2$ do not occur in either of the paths
$p_1$ or $p_2$: in order to use either in $p_i$, the path would first have to
return to $x_0$ and then leave again on $e_i$, contradicting minimality.
By Lemma~\ref{compacts}, there are (contractive) injective operators $S_1$
and $S_2$ on $\B(\H)$ so that $\rho(c_1) S_1$ and $\rho(c_2) S_2$ generate the
compact operators as an algebra.
Replace $K_1$ and $K_2$ with $K_1 S_1$ and $K_2 S_2$ and redefine $\rho$ with
the new values.
We still have the \wot-continuous extension to $\fL_G$.
But in addition, we now have that
\[
 \rho( P_{x_0} \fL_G P_{x_0}) = Q_{x_0} \rho(\fL_G) Q_{x_0}
\]
contains generators for $\fK(\H_{x_0})$.

Let us show that $\rho$ is irreducible.
The \wot-closure of the range of $\rho$ contains $\B(\H_{x_0})$.
If $x$ is any other vertex of $G$, select a path $p$ from $x_0$ to $x$
and a path $q$ from $x$ to $x_0$.
Then $\rho(p)$ is a compact operator which is injective on $\H_{x_0}$
and maps this onto a dense subspace of $\H_x$.
Likewise $\rho(q)$ is a compact operator which is injective on $\H_x$
and maps it onto a dense subspace of $\H_{x_0}$.
Therefore $\rho( P_x \fL_G P_{x_0})$ contains $Q_x \rho(p) \fK(\H_{x_0})$
which is norm-dense in $Q_k \fK  Q_{x_0}$, and thus \wot-dense in
$Q_x \B(\K) Q_{x_0}$.
Similarly $\rho( P_{x_0} \fL_G P_x)$ is \wot-dense in $Q_{x_0} \B(\K) Q_x$.
Whence $\rho( P_x \fL_G P_x)$ is \wot-dense in $Q_x \B(\K) Q_x$ as well.
It follows that $\rho$ is irreducible.

Next we verify that $\rho$ is faithful.
We wish first to show that the restriction to $Q_{x_0}\K$ is
faithful on $P_{x_0} \fL_G P_{x_0}$.
Enumerate all minimal cycles with source $x_0$ as $\C = \{ c_1, c_2, \dots \}$
beginning with the two cycles $c_1$ and $c_2$ used above.
Now $x_0 \fgeeplus x_0$ is isomorphic to $\bF^+(G_{x_0})$ where $G_{x_0}$ is
the graph with a single vertex $x_0$ and loop edges given by $\C$.
In particular, $P_{x_0} \fL_G P_{x_0}$ is isomorphic to $\fL_{G_{x_0}}$, which
is in turn isomorphic to $\fL_n$ where $n$ is the cardinality of $\C$.
Following the proof of Theorem~\ref{Ln_irred}, as we have the conditions of
dense range and generation of the compact operators, it suffices to verify
that $\sum_{i\ge1} \|\rho(c_i)\| = r < 1$.
There are a number of ways to see this.
One is this: the calculation of $\|V\|$ in the proof of Theorem~\ref{Ln_irred}
bounds the sum of the upper bounds of $\|\rho(w)\|$ over all paths by
$(1-r)^{-1}$.  The corresponding calculation using the same norm estimates for
$\rho(c_i)$ must yield a smaller sum, as the collection of paths built from
$\C$ is a proper subset.  But if $\sum_{i\ge1} \|\rho(c_i)\| > r$, this
estimate would be strictly greater.
We can conclude that $\rho$ is faithful on $P_{x_0} \fL_G P_{x_0}$.

Suppose that $A$ is a non-zero element in $\ker \rho$.
Then for some vertices $x$ and $y$, $P_y A P_x \ne 0$.
Choose a path $p$ from $x_0$ to $x$ and a path $q$ from $y$ to $x_0$.
Then $L_q A L_p$ is a non-zero element of $P_{x_0} \fL_G P_{x_0}$ in $\ker
\rho$, contrary to fact.
Hence $\rho$ is faithful.
\bx

\begin{rem}\label{transpose}\textit{The transpose graph.}
We need to make an observation about representations of the transposed graph
$G^t$ which has the same vertices and edges as $G$, but the edges have the
range and source maps reversed.
The semigroupoids $\bF^+(G)$ and $\bF^+(G^t)$ are anti-isomorphic but
generally this does not extend to the algebras $\fL_G$ and $\fL_{G^t}$.
Kribs and Power \cite{KrP} show that $\fL_{G^t}$ is unitarily equivalent to
$\fR_G$, the algebra generated by the right regular representation.
However the fact that we are in a strictly contractive situation makes things
easier.

Observe that the diagonal algebra
$\D = \textsc{wot--}\!\!\spn\{P_x : x \in \V(G) \}$ is
mapped isometrically onto the diagonal
$\textsc{wot--}\!\!\spn\{ \rho(x) : x \in \V(G) \}$.
Consider the \wot-closed ideal $\fL_G^0$ generated by $\{L_e : \in \E(G) \}$.
Note that $\fL_G = \D + \fL_G^0$ as a topological direct sum.
Let $A(\bF^+(G))$ denote the (separable) algebra of absolutely convergent sums
of  $\{ L_w : w \in \bF^+(G) \}$ with norm
$\| \sum a_w L_w \|_1 = \sum_{w \in \bF^+(G)} |a_w|$.
The restriction of $\rho$ to $\fL_G^0$ maps factors the ideal through
$A(\bF^+(G))$.  Indeed, set $r_e = \|K_e\|$.
The map determined by $\sigma(L_e) = r_e L_e$ maps $\fL_G^0$ into
$A(\bF^+(G))$ because Fourier coefficients are bounded and
\begin{align*}
 \sum_{k \ge 1} \sum_{|w|=k} \| \sigma(L_w)\| &\le
 \sum_{k \ge 1} \sum_{|w|=k} r_{e_{i_k}} \dots r_{e_{i_1}} \\&\le
 \sum_{k \ge 1} \big( \sum_{e \in \E(G)} r_e \big)^k \le \frac r{1-r}
\end{align*}
where we write $w = e_{i_k} \dots e_{i_1}$.
Then define a representation $\tau$ of $A(\bF^+(G))$ by setting
$\tau(L_e) = \| \rho(L_e) \|^{-1} \rho(L_e)$ and $\tau(x) = \rho(x)$.
Clearly this determines a contractive representation of $A(\bF^+(G))$.
Moreover the factorization $\rho = \tau \sigma$ is evident.

The representation $\tau$ inherits certain properties from $\rho$.
For example, since the \wot-closure of the range of $\tau$ equals
the \wot-closure of the range of $\rho(\fL_G)$, it will be irreducible if
and only if $\rho$ is.   Moreover in the construction used here, it will
be faithful if $\rho$ is faithful.  This essentially comes down to
a calculation as in the proof of Theorem~\ref{Ln_irred} that shows that the
intertwining operator $V$ is injective.

Now observe that $A(\bF^+(G^t))$ is anti-isomorphic to
$A(\bF^+(G))$ in the natural way by identifying a path with the reverse.
The representation
\[
 \tau^t \big( \sum_{w \in \bF^+(G^t)}a_w L_w \big) =
 \sum_{w \in \bF^+(G^t)} a_w \tau(L_{w^t})^*
\]
is continuous, and has range equal to $\tau( A(\bF^+(G)) )^*$.
In particular, if $\tau$ is injective, then so is $\tau^t$; and
if $\tau$ is irreducible, so is $\tau^t$.

Now consider the representation $\rho^t = \tau^t \sigma^t$ of $\fL_{G^t}$.
This is the representation of $\fL_{G^t}$ which satisfies
$\rho^t(L_w) = \rho(L_{w^t})^*$.
This extends to a \wot-continuous completely contractive representation of
$\fL_{G^t}$ by Corollary~\ref{dilationcor}.
Now $\sigma^t$ maps $\fL_{G^t}$ injectively into a dense subalgebra of
$A(\bF^+(G^t))$.
If $\rho$ is faithful, then so is $\tau$; whence $\tau^t$ and $\rho^t$ are
also faithful.  Likewise $\rho^t$ will be irreducible if $\rho$ is
irreducible.  For the next application, we observe that if the range of
$\rho$ is a nest representation, then so is $\rho^t$ for the complementary
nest.
\end{rem}

We can now establish a complete characterization of when a graph
algebra has a faithful nest representation.

\begin{defn}\label{transitive_quotient}
If $G$ is a countable directed graph, the \textit{transitive quotient} of $G$
is the graph obtained by replacing each transitive component with a vertex
$\tau$, each edge mapping from one component to another becomes an edge between
the corresponding vertices of the quotient graph, and the number of loop edges
at a vertex $\tau$ will be the number of primitive cycles through any fixed
vertex in the corresponding component.
\end{defn}

Observe that there is no loop on a vertex $\tau$ in the transitive quotient
only when the transitive component is trivial.  There is a single loop edge
precisely when the component is a cycle.  When a component is strongly
transitive, there will be countably many loop edges in the quotient.

\begin{thm}\label{faithful_nest_rep}
 A directed graph $G$ has a faithful nest representation if and only if
\begin{enumerate}
 \item The transitive quotient of $G$ is totally ordered.
 \item No transitive component of $G$ is a cycle.
 \item The trivial components of $G$ form an interval that is order
   equivalent to a subset of $\bZ$ in the total order on $G$,
   and the restriction of the graph to this set consists of a single
   path with one edge from each vertex to the next in the order.
\end{enumerate}
\end{thm}

\Prf First observe that these conditions are necessary for the existence of
a faithful nest representation $\rho$ into $\T(\N)$.
It is evident that each transitive component must map onto an interval of the
nest, and that this restriction must be faithful.  So cycle components are
eliminated as they yield infinite dimensional abelian algebras.
Moreover there must be edges from one component to another if and only if the
image interval dominates the other interval in the order on the nest.
Hence the transitive quotient is totally ordered.

If $x$ and $y$ represent trivial components, then $\rho(x)$ and $\rho(y)$ are
one dimensional.
Suppose that $\rho(x) \prec \rho(y)$ in the order on the nest.
Then $\rho(x) \T(\N) \rho(y) = \rho( P_x \fL_G P_y)$ is one dimensional.
Hence $ P_x \fL_G P_y$ is one dimensional.
Thus there is a unique path from $y$ to $x$ in $G$.
So this path cannot pass through any strongly transitive component.
So the path passes only through trivial components.
Moreover, there can only be one edge from a trivial component leading to
another trivial component (for otherwise from the total order, we would deduce
the existence of two paths from one vertex to another).
In particular, the order structure on the trivial components must coincide
with a subset of $\bZ$.

For the converse, first suppose that every transitive component is strongly
transitive and that the transitive quotient is totally ordered.
We will make a modification in the construction of an irreducible
representation as follows.
Add some extra edges to $G$ to form a (strongly) transitive graph $H$.
Use Lemma~\ref{op_range} to define an injective compact operator
$K = \begin{bmatrix}K_1 & K_2 & \dots\end{bmatrix}$ so that each $K_i$ has
dense range and $\sum_{i\ge1} \|K_i\| = r < 1$.
As in the proof of Theorem~\ref{faithful_irred}, we can now define a
faithful irreducible representation $\phi$ of $\fL_H$.
This is not affected if certain of the operators $K_i$ are replaced by
$K_i S_i$ if $S_i$ is an invertible contraction.
So using Lemma~\ref{compacts} as in the proof of
Theorem~\ref{faithful_irred}, we can modify some of these operators so that
the restriction to each transitive component of $G$ is \wot-dense in the
algebra of operators on its range. Now restrict the representation to $\fL_G$.

Clearly the restriction is faithful.
We claim that $\phi(\fL_G)$ is \wot-dense in a nest algebra
$\T(\N)$ which has an infinite dimensional atom for each transitive component
of $G$ and is ordered as the transitive quotient.
Indeed, for each transitive component $\tau$, let
\[
  P_\tau = \sotsum_{x \in \V(\tau)} \phi(x) .
\]
Then by the construction above, we have arranged that $\phi(\tau)$ is
\wot-dense in $\B(P_\tau\H)$.
Form the nest $\N$ obtained by ordering the subspaces $P_\tau\H$ by the total
order on the transitive quotient of $G$.
Observe that the range of $\phi$ is \wot-dense in the diagonal algebra
\[
 \T(\N) \cap \T(\N)^* = \sum_{\tau}\! \oplus \B(P_\tau\H).
\]
Finally, from the fact that the induced order on the transitive quotient is
total, there will be a path from each component to every other lower in the
total order.
The image of this path yields a non-zero operator between the corresponding
atoms of the nest.  Together with the diagonal, it is straightforward to
check that the algebra generated is \wot-dense in $\T(\N)$.

To deal with an interval of trivial components, proceed as follows.
Let $\V_t$ denote the vertices of $G$ which are trivial components.
Let $\V_+$ denote the vertices of $G$ which are greater than $\V_t$
in the total order, and $\V_-$ be those which are less than $\V_t$.  Let
$G_-$, $G_t$ and $G_+$ be the restriction of $G$ to $\V_-$, $\V_t$ and $\V_+$
respectively. In addition, the restriction of $G$ to $\V_t \cup V_-$
and $\V_t \cup V_+$ will be denoted as $G_{-,t}$ and $G_{t,+}$
respectively.

Form a Hilbert space $\H =  \sumoplus_{x \in \V(G)} \H_x$
where $\H_x$ is infinite dimensional for $x \in \V_- \cup \V_+$
and $\H_x = \bC h_x$ is one dimensional for $x \in \V_t$.
The representation will have $Q_x := \phi(x) = P_{\H_x}$ for $x \in \V(G)$.
As before, identify each $\H_x$ for $x \in \V_- \cup \V_+$ with $\H$ via a
unitary $J_x$ of $\H$ onto $\H_x$.
Again start with an injective compact operator $K = [ K_e ]_{e \in \E(G)}$
such that each component $K_e$ has dense range and
\[
 \|K\| \le \sum_{e \in \E(G)} \|K_e\| \le r < 1 .
\]

We define the representation $\phi$ in several stages. First, we define $\phi$
on $G_{-,t}$.
Use the construction above for the case of all strongly transitive components
to define a faithful nest representation $\phi^-$ of $\fL(G^-)$.
This involved some modification of the $K_e$'s.

Enumerate the vertices $x_j$ in $\V_t$ in their given order by a subset
$J \subset \bZ$; and let $e_j$ be the unique edge in $G_t$ with
$s(e_j) = x_j$ and $r(e_j) = x_{j-1}$ for all $j\in J$ except the least
element (if there is one).
Define
\[
 \phi(e_j) = r_j h_{j-1} h_j^*
\]
where $r_j = \| K_{e_j}\|$.
Here the role of the operator $K_{e_j}$ is merely to control the norm
of $\phi(e_j)$.

Set
\[
 \E_- = \{ e \in \E(G) : s(e) \in \V_t \AND r(e) \in \V_- \}
\]
and dually
\[
 \E_+ = \{ e \in \E(G) : s(e) \in \V_+ \AND r(e) \in \V_t \} .
\]
Fix a unit vector $u_0 \in \H$.
Define a vector $k_e = J_{r(e)}K_e u_0$ for each $e \in \E_-$;
and let $\phi(e) = k_e h_x^*$. This defines $\phi$ on $G_{-,t}$.

For $G_{t,+}$, we will instead use the same construction for the transposed graph
$G_{t,+}^t$, and take its transpose as in Remark~\ref{transpose}. Finally,
if $e$ is an edge with $s(e) \in \V_+$ and $r(e) \in \V_-$, define
$\phi(e) = J_{r(e)} K_e J_{s(e)}^*$. This defines $\phi(e)$ for all edges of $G$.
Since
\[
 \sum_{e \in \E(G)} \|\phi(e)\| \le
 \sum_{e \in \E(G)} \|K_e\| \le r < 1 ,
\]
this defines a \wot-continuous completely contractive representation of
$\fL_G$ by Corollary~\ref{dilationcor}.

The nest $\N$ is built from these atoms using the total order
on the transitive quotient as before.

Observe that this is a nest representation.
The restriction of $\phi$ to $G_-$ is a faithful nest representation of
$\fL_{G_-}$ by the earlier construction.
The restriction of $\phi$ to $G_+$ is a faithful nest representation
of $\fL_{G_+}$ because of Remark~\ref{transpose}.
The restriction of $\phi$ to $G_t$ maps $\fL_{G_t}$ into a
\wot-dense subalgebra of the upper triangular operators on the basis
$\{h_{x_j} : j \in J \}$.
This is also faithful because there are no extra edges involved.
Moreover there are edges mapping $G^+$ into $G^t$ and
from $G^t$ into $G^-$ so that there is a path from every vertex in
$\V_+$ to any vertex in $\V_t$ and from every vertex in $\V_t$
to every vertex in $\V_-$.
Consequently, the fact that the edges are sent to non-zero operators is
sufficient to see that the algebra generated is \wot-dense in the nest
algebra.

It remains to verify that the map $\phi$ is faithful.
Clearly, $\ker \phi$ is generated by elements of the form $A = P_y A P_x$
where $x, y \in \E (G)$.
So it suffices to show that $\phi$ does not annihilate elements of this form.
Since the restrictions to the three subgraphs are faithful, we need only
consider three possibilities.

Consider a non-zero element $A \in \fL_G$ of the form $A = P_y A P_x$ for
$x \in \V_t$ and $y \in \V_-$.
Let $A \sim \sum_{w=ywx} a_{w} L_w$ be its Fourier expansion.
Each $w$ in the above sum factors uniquely as $w=w'ew''$, with
$w' \in \bF^{+}(G_-)$, $e \in \E_{-}$ and $w'' \in \bF^{+}(G_t)$.
Let $k$ be the length of the shortest path $w'$ appearing in these
factorizations, say $\tilde{w}_0 = w_0 e_0 w_0''$.

Define a set of paths $W$ to be
\[
  \{ we \in \bF^{+}(G) : w \in \bF^{+}(G_-),\, |w|=k,\ r(w)=y,\,
             e \in \E(G_-)\!\cup\! \E_- \}.
\]
Then $A$ factors as:
\[
  A = \sum_{we \in W} L_{we} A_{we}
    = \big[ L_{we} \big]_{we \in W} \big[ A_{we} \big]_{we \in W}^t
\]
where $A_{we} = L_{we}^* A$.
This is a bounded factorization in which $\big[ L_{we} \big]_{we \in W}$ is a
partial isometry and $\| \big[ A_{we} \big]_{we \in W}^t \| = \|A\|$.
\vspace{.5ex}
The term $L_{w_0 e_0} A_{w_0 e_0}$ is non-zero.
We may suppose that $A_{w_0 e_0} = P_{s(e_0)} A_{w_0 e_0} P_x$, and this
is a non-zero element of $\fL_{G_t}$.

But then, $\phi(A_{w_0 e_0}) h_x$ is a non-zero multiple of $h_{s(e_0)}$,
which is send by $\phi(L_{e_0})$ to a non-zero multiple of $k_{e_0}$, and
then by $\phi(L_{w_0})$ to a non-zero vector in the range of
$J_y K_{w_0 e_0}$.  On the other hand,
\[
 \phi(A - L_{w_0 e_0} A_{w_0 e_0})  =
 \big[\phi(L_{we}) \big]_{we \in W \bsl \{w_0e_0\}}
 \big[ \phi(A_{we}) \big]_{we \in W \bsl \{w_0e_0\}}^t
\]
lies in the range of $J_y \big[ K_{we} \big]_{we \in W \bsl \{w_0e_0\}}$.
The intersection is $\{0\}$ since $J_y \big[ K_{we} \big]_{we \in W}$ is
injective.
Therefore $\phi(A) \ne 0$.
So $\phi$ is injective on the graph restricted to $\V_- \cup \V_t$.

A similar argument shows that $\phi$ is injective on the graph restricted to
$\V_- \cup \V_t$.
Indeed, taking adjoints puts us in the same case as above for the transpose
algebra. This is why we used the transpose of a representation of
$G_{t,+}^t$. To return to $G_{t,+}$, we use Remark~\ref{transpose}.

Finally we consider non-zero elements $A = P_y A P_x$ with $x \in \V_+$ and
$y \in \V_-$.  The argument is much the same.
Again there is a least integer $k$ so that $A$ has a non-zero Fourier
coefficient $a_{w_0e_0w_0''}$ where $w_0 \in \bF^+(G_-)$,
$|w_0|=k$ and $s(e_0) \in \V_+\cup\V_t$.
Again define a set of paths $W$ as
\[
 \{ we \in \bF^{+}(G) : w \in \bF^{+}(G_-),\, |w|=k,\ r(w)=y,\,
             s(e) \in \V_+\cup\V_t \}.
\]
Consider once again factorizations of the form
\[
  A = \sum_{we \in W} L_{we} A_{we}
    = \big[ L_{we} \big]_{we \in W} \big[ A_{we} \big]_{we \in W}^t
\]
where $A_{we} = L_{we}^* A$.
This is a bounded factorization in which $\big[ L_{we} \big]_{we \in W}$ is a
partial isometry and $\| \big[ A_{we} \big]_{we \in W}^t \| = \|A\|$.
\vspace{.5ex}
The term $L_{w_0 e_0} A_{w_0 e_0}$ is non-zero.
We may suppose that $A_{w_0 e_0} = P_{s(e_0)} A_{w_0 e_0} P_x$, and this
is a non-zero element of $\fL_{G}$.

But then, $\phi(A_{w_0 e_0})$ has non-zero range in $\H_{s(e_0)}$,
which is send by the injective operator $\phi(L_{w_0 e_0})$ to a non-zero
vector in the range of $J_y K_{w_0 e_0}$.  On the other hand,
\[
 \phi(A - L_{w_0 e_0} A_{w_0 e_0})  =
 \big[\phi(L_{we}) \big]_{we \in W \bsl \{w_0e_0\}}
 \big[ \phi(A_{we}) \big]_{we \in W \bsl \{w_0e_0\}}^t
\]
lies in the range of $J_y \big[ K_{we} \big]_{we \in W \bsl \{w_0e_0\}}$.
The intersection is $\{0\}$ since $J_y \big[ K_{we} \big]_{we \in W}$ is
injective.
Therefore $\phi(A) \ne 0$.

This completes all of the cases and establishes that $\phi$ is faithful.
\bx

\begin{eg}\label{Cantor}
The nest used in the this theorem need not have a simple order.
For example, consider a graph $G$ with vertices indexed by the rationals $\bQ$.
Put two loop edges at each vertex and an edge from $x_s$ to $x_r$ if $r<s$.
Then the theorem establishes a faithful representation into the Cantor nest
with infinite dimensional atoms.
This nest is similar to a nest with a large continuous part.
Thus there are \wot-continuous completely bounded faithful nest
representations into nests which are not atomic.

We do not know which graphs have a faithful representation into a continuous
nest algebra.
\end{eg}

The final result determines when a free semigroupoid algebra $\flgee$
admits a nest representation for a multiplicity-free nest which has
basis of atoms  ordered by $\bN$.
It is remarkable that in spite of the large radical in the nest algebra,
if $G$ has no sinks, then such a free-semigroupoid algebra has to be
semisimple.
This result can be modified to yield a characterization for nests of order
type $-\bN$ and $\bZ$ as well, but the details are omitted.

\begin{thm} \label{Nrepn}
Let $G$ be a countable directed graph and let $\N$ be the multiplicity-free
nest which is ordered as $\bN \cup \{+\infty\}$. The free semigroupoid algebra
$\flgee$ admits a faithful nest representation into $\Alg \N$ if and only
if one of the following conditions is satisfied:
\begin{enumerate}

\item $G$ is strongly transitive and every vertex $x \in \V(G)$ supports a
 loop edge.

\item $\V(G) = \{ x_i : i \ge 1 \}$ and $\E(G) = \{ e_i : i \ge 1 \}$
 where $s(e_i) = x_{i+1}$ and $r(e_i) = x_i$ for all $i \ge 1$.

\item $\V(G) = \V_0 \cup \{ x_i : 1 \le i \le n \}$.  The restriction of $G$
 to $\V_0$ is a strongly transitive graph $G_0$ such that every vertex
 $x \in \V_0$ supports a loop edge. The restriction of $G$ to
 $\{ x_i : 1 \le i \le n \}$ consists of $\{ e_i : 1 \le i \le n-1 \}$ where
 where $s(e_i) = x_{i+1}$ and $r(e_i) = x_i$ for all $1 \le i \le n-1$.
 Plus there is an edge $e_n$ with $r(e_n) = x_n$ and $s(e_n) \in \V_0$.
 There may, in addition, be other edges with source in $\V_0$ and range in
 $\{ x_i : 1 \le i \le n \}$.

\end{enumerate}
\end{thm}

\Prf The nest $\N$ is unique up to unitary equivalence, so let $\H = \ell^2$
with orthonormal basis $\{ h_i : i \ge 1 \}$, and set
\[
 N_k = \spn\{ h_i : 1 \le i \le k \} \qfor k \ge 0 \qand N_\infty = \H .
\]
The nest is
\[
 \N = \big\{ N_k : k \in \bN_0 \cup \{\infty\} \big\} .
\]
The atoms of the nest are the projections
$P_i = h_i h_i^* = P_{N_i} - P_{N_{i-1}}$ for $i \ge 1$.

Assume that there exists a faithful nest representation $\phi$ mapping
$\flgee$ onto a \wot-dense subalgebra of $\Alg \N$. Then vertices are
sent to diagonal projections $Q_x := \phi(x)$ for $x \in \V(G)$.

As in the proof of Theorem~\ref{upper} about upper triangular nest
representations, we find that if $x$ supports no loop edge, then $Q_x$
will have rank at most one (thus equal to one since $\phi$ is faithful).
But then this forces $x \fgeeplus x$ to be one dimensional, which implies that
there are no cycles through $x$.  Thus each transitive component that consists
of more than one point has a loop edge at each vertex.

In the case that the transitive component of $x$ is not the trivial graph,
there are cycles through $x$.  Therefore $x \fgeeplus x$ is infinite
and $P_x \fL_G P_x$ is infinite dimensional.
Therefore $\phi(P_x)$ has infinite rank.
This is a diagonal projection, so its support must be unbounded.
Therefore it follows that $\phi(P_y) \Alg \N \phi(P_x)$ is non-zero for every
edge $y$.  Consequently there is a path in $G$ from $x$ to $y$.
So the infinite rank points form a single transitive component.

There is a difference from the situation of Theorem~\ref{upper} regarding a
transitive component consisting of exactly one point.  We claim that in this
instance, there are at least two loop edges or none at all.  In other words,
it cannot be a cycle consisting of a single point and single loop edge.
Indeed, if $x$ is itself a transitive component with a single loop edge $u$,
then $x \fgeeplus x$ consists of powers of $u$ and $P_x \fL_G P_x$ is
isomorphic to $H^\infty$.  This is an infinite dimensional abelian algebra.
As in the previous paragraph, $\phi(x)$ has infinite rank.
But then $\phi(x) \Alg \N \phi(x)$ is not abelian,
so it cannot be contained in the range of $\phi( P_x \fL_G P_x)$.
This establishes our claim.

So if $G$ contains no rank one vertices, it consists of a single strongly
transitive component in which each vertex supports a loop edge.
This is case (1).

In the event that $\phi(x)$ has rank one, it must equal some atom
$P_n$ of the nest algebra.
Hence $\Alg \N \phi(x)$ is finite dimensional.
So $\fL_G P_x$ is finite dimensional.
This implies that $x$ has no paths into any cycles, and thus the path
component beginning at $x$ is finite with no loops.
Thus no atom $P_i$ for $i<n$ is dominated by $\phi(y)$ for any $y$
in a strongly transitive component.  So all correspond to rank one vertices.
The vertices $x_i$ which have rank one are linearly ordered in the
transitive relation determined by $\Alg \N |_{\sum \phi(P_{x_i})\H}$
because exactly one of $\phi(P_{x_j}) \Alg \N \phi(P_{x_i})$ or
$\phi(P_{x_i}) \Alg \N \phi(P_{x_j})$ is non-zero for $i \ne j$.

Say that an initial segment $P_1, \dots, P_n$ of $\N$ are atoms of the nest
(corresponding to unit vectors $h_1,\dots,h_n$)
which are the image of rank one vertices $x_1,\dots, x_n$.
For $2 \le k \le n$, There is an element $a_k \in \fL_G$ such that
$\lip \phi(a_k) h_k, h_{k-1} \rip = 1$.
We replace $a_k$ by $P_{x_{k-1}} a P_{x_k}$ so that
$\phi(a_k) = h_{k-1} h_k^*$.
It follows that there is an edge $e_k$ in the graph from $x_k$ to
$x_{k-1}$, for any other path will factor through other vertices
and so would not be supported on the $(k-1,k)$ entry.

Moreover it follows that there are no other edges between rank one vertices.
To see this, suppose that there is an edge $e$ from $x_k$ to $x_j$
where $j \le k-2$.
Then arguing as above, we see that $\phi(e) = \alpha h_j h_k^*$ for
some $\alpha \in \bC$.
Now $a = \alpha a_{j+1}\dots a_k$ is a word of order at least two.
Note that $\phi(a - e) = 0$ and $a \ne e$, contradicting the
faithfulness of $\phi$.
It follows that the graph on $x_1,\dots,x_n$ is a simple directed chain.

If $G$ consists entirely of rank one vertices, then $G$ is the graph
corresponding to case (2).

If there is a non-empty finite collection $x_1,\dots,x_k$ of rank one
vertices, then there is a single strongly transitive component which is mapped
into a dense subalgebra of $P_{N_k}^\perp \Alg \N$.  A particular vertex $y_1$
must satisfy $\phi(y) > P_{n+1}$.  Arguing as above, there must be an edge in
$G$ from $y$ to $x_n$.  But there may also be other edges from the transitive
component to the rank one vertices. This is case (3).

Thus we have shown that one of the three possible types of graph above
are necessary for a nest representation onto a dense subalgebra of $\Alg \N$.

For the converse, we deal separately with the three cases.

Case (2) is simplest.  Just define $\phi(x_i) = h_i h_i^*$ and
$\phi(e_i) = \frac12 h_i h_{i+1}^*$.  This extends to a \wot-continuous
representation of $\fL_G$ by Corollary~\ref{dilationcor}.
The range is evidently upper triangular and contains the matrix units.
Hence it is \wot-dense in $\Alg\N$.

Now suppose that (1) holds.
We use the notation and terminology from the proofs of
Theorem~\ref{irreducible} and Theorem \ref{upper}.

As in the proof of Theorem~\ref{upper}, we select a loop edge $f_x$ for each
vertex $x \in \V(G)$.  Let $\L = \{ f_x : x \in \V(G) \}$.
As we have observed in Theorem~\ref{upper} and Corollary \ref{cor:Fourier},
there exists a countable family $\{ \psi_{w_n,\lambda_n} \}_{n\ge1}$ of
upper triangular representations of $\flgee$ which separates the points in
$\flgee$.
These paths $w_n$ do not involve the edges $f_x$.
A moments reflection shows that the scalars $\lambda_{n,i} \in \bT$
used in these representations may be chosen to be distinct.

Choose an infinite path in $G$ of the form $w = e_{j_1} e_{j_2} e_{j_3} \dots$
such that $x_0=r(e_{j_1})$ and $x_k = r(e_{j_k}) = s(e_{j_{k-1}})$ for $k \ge
1$. We require that no edge of the form $f_x$ is used, and that each path $w_n$
should occur separately somewhere within the path $w$.
Assign distinct scalars $\lambda_k \in \bT$ for $k\ge0$ such that on each
segment corresponding to one of the words $w_n$, we use the scalars
$\lambda_{n,i}$ selected above for the separating family of representations.
Write $\lambda = (\lambda_k)_{k\ge0} \in \bT^{\bN_0}$.

Define a representation $\psi_{w,\lambda}$ of $\fgeeplus$ into $\T(\N)$ as
follows:
{\allowdisplaybreaks
\begin{alignat*}{2}
 \psi_{w,\lambda}(x)   &= \sum_{x_k=x} h_k h_k^*
       &\FOR& x \in \V(G) \\
 \psi_{w,\lambda}(f_x) &= \frac 1 2 \sum_{x_k=x} \lambda_k h_k h_k^*
       &\FOR& x \in \V(G) \\
 \psi_{w,\lambda}(e)   &= \frac 1 2 \sum_{e_{j_k}=e}h_{k-1} h_k^* \quad
       &\FOR& e \in \E(G) \bsl \L
\end{alignat*}
}  

Arguing as in Theorem~\ref{upper}, the row operator
\[
  E = \begin{bmatrix}\psi_{w,\lambda}(e)\end{bmatrix}_{e \in \E(G)}
\]
satisfies $EE^* \le \frac12 I$.
So by Corollary~\ref{dilationcor}, this representation extends to a
\wot-continuous completely contractive representation of $\fL_G$.

Observe that the image of each vertex is diagonal and the image of each edge
is upper triangular, so the range of this representation is contained in
$\T(\N)$.
The algebra generated by $\{ \psi_{w,\lambda}(f_x) : x \in \V(G) \}$ contains
a family of diagonal operators with distinct non-zero entries in each
coordinate.  Consequently, the \wot-closed algebra that they generate is the
full diagonal algebra $\D$.
This together with the rank one partial isometries
\[
 h_{k-1}h_k^* = \psi_{w,\lambda}(e_{j_k})
\]
generates $\T(\N)$ as a \wot-closed algebra.  Therefore this is a nest
representation.

Finally observe that for each $n$, there is an interval of the nest
corresponding to the representation $\psi_{w_n,\lambda_n}$.
The compression to this interval is equivalent to it.
As the direct sum of these representations is faithful, it follows that
$\psi_{w,\lambda}$ is also faithful.

It remains to deal with case (3).  There is a finite initial segment
of the graph with vertices $x_1,\dots,x_n$ and edges $e_1,\dots,e_n$
where $r(e_i) = x_i$ and $s(e_i) = x_{i+1}$ for $1 \le i < n$,
$r(e_n) = x_n$ and $s(e_n) = y_0$ belongs to the strongly transitive
component.  Let $G_0$ be the graph $G$ restricted to $\V_0$.

There is a faithful upper triangular nest representation of $\fL_{G_0}$
determined by the previous paragraphs.
We will need to specify some special features of this representation.
Let $\L = \{ f_x : x \in \V(G_0) \}$ be a choice of a loop edge at each vertex
of $G_0$.
For each word $v_s$ built from edges in $\E(G_0) \bsl \L$, choose a
countable family of points
$\lambda_{v_s,k} = ( \lambda_{v_s,k,i} )_{i=0}^{|v_s|}$ dense in
$\bT^{|v_s|+1}$.
We may also ensure that all $\lambda_{v_s,k,i}$ are distinct.
For later use, we split the collection of all these pairs
$(v_s, \lambda_{v_s,k})$ into countably many families $U_j$ with the same
properties.

Choose an infinite word $w$ with $r(w) = y_0$ which contains each word
$v_s \in \bF^+(G_0) \bsl \L$ infinitely often, and denote the $k$th
occurrence by $v_{s,k}$.
Also choose $\lambda = (\lambda_j)_{j\ge1} \in \bT^\bN$ so that the
coordinates $\lambda_j$ are distinct, and coincide with $\lambda_{v_s,k,i}$
at the occurrence of the word $v_{s,k}$ in $w$.
Denote the pairs $v_{s,k}$ as $u_i$ in order of occurrence.
We may assume that they are sufficiently widely spaced that
there are more than $4i$ edges between $u_i$ and $u_{i-1}$.
We obtain an upper triangular nest representation $\psi_{w,\lambda}$ as above.
Denote the basis on which this acts as $\{h_{n+j} : j \ge 1 \}$.

The plan is to extend this representation to all of $G$.
Part of the extension is natural.  Set
{\allowdisplaybreaks
\begin{alignat*}{2}
 \psi(p) &= \psi_{w,\lambda}(p) &\FOR& p \in \bF(G_0) \\
 \psi(x_i)   &= h_i h_i^*
       &\FOR&1 \le i \le n \\
 \psi(e_i) &= \frac12 h_i h_{i+1}^* \quad
       &\FOR& 1 \le i \le n-1
\end{alignat*}
}  
We also need to define $\psi$ on edges in $G$ with source in $\V_0$ and range
in the initial segment including $e_n$.
Let $\{ g_j : j \in J \}$, $J \subset \bN_0$, be an enumeration of these
edges including $g_0=e_n$.
Let $r(g_j) = x_{i_j}$ and  $s(g_j) = y_j \in \V_0$,
We will define vectors $k_j = \psi(y_j) k_j$ and set
\[
 \psi(g_j) = h_{i_j} k_j^* \qfor j \in J .
\]
We need to define the vectors $k_j$ in such a way that the representation
is faithful.
At this point, we specify only that $\alpha := \ip{k_0,h_{n+1}} \ne 0$.

Before dealing with that issue, let us observe that whatever the choice
of vectors $k_j$, this representation will be a nest representation.
Indeed, the compression to $\spn \{ h_{n+k} : k \ge 1 \}$ is a faithful
map of $\fL_{G_0}$ into a \wot-dense subalgebra of the upper triangular
nest  algebra. The restriction to $\spn \{ h_1 : 1 \le i \le n \}$ is a
faithful map of the interval onto $\T_n$.
Moreover $\rho(e_n) h_{n+1} h_{n+1}^* = \frac{\alpha}{2} h_n h_{n+1}^*$
provides an element in the range which now generates the full nest algebra.
The addition of more edges only adds to this already dense subalgebra.

For each $j \in J$, set
\[
 V_j = \{ i \ge 1 :  u_i \in U_j,\ r(u_i) = y_j \}.
\]
The vertex $r(u_i)$ is identified with a basis vector $h_{r_i}$ in
the representation $\psi$, and $s(u_i)$ is identified with $h_{s_i}$
where $s_i = r_i + |u_i|-1$.
Define
\[
 k_0 = \frac12 h_{n+1} + \sum_{i \in V_0} 2^{-i-1} h_{r_i}
\]
and
\[
 k_j = \sum_{i \in V_j} 2^{-i-1} h_{r_i} \qfor j \in J \bsl \{0\}.
\]
Observe that $k_j = \psi(y_j) k_j$ for each $j \in J$.

Again we can verify that the row operator
$E = \big[ \psi(e) \big]_{e \in \E(G)}$ satisfies $EE^* \le \frac12 I$.
Therefore by Corollary~\ref{dilationcor}, $\psi$ extends to a completely
contractive \wot-continuous representation of $\fL_G$.

The last step is to verify that $\psi$ is faithful.
We have already shown that it is faithful on $\fL_{G_0}$ and on the initial
segment.
So it suffices to consider elements with domain in $G_0$ and range in the tail.
This is the ideal $\J$ generated by the $\{ g_j : j \in J \}$.

To this end, suppose that $A \in \J$ satisfies $\psi(A)=0$.
This is equivalent to $h_i h_i^* \rho(A) = \rho(P_{x_i}A) = 0$ for
$1 \le i \le n$.
So it suffices to consider $A = \sum_{j \in J} q_j g_j A_j$
where $q_j$ is the unique path in $G$ from $r(g_j)$ to $x_i$
(if $i_j \ge i$, else it is $0$), and $A_j = P_{y_j} A_j$
lies in $P_{y_j} \fL_{G_0}$.
The range of $\psi(A)$ is contained in $\bC h_i$.
So $\psi(A) = 0$ if and only if
\[
 0 = \psi(A)^* h_i = \sum_{j \in J} \psi(A_j)^* k_j .
\]

We can write the Fourier series of each term as
$A_j \sim \sum_{r(v) = y_j} a_{v,j} L_v$.
Also note that in $\fL_G$, we have $\|A\|^2 = \| \sum_{j \in J} A_j^* A_j \|$.
As in the previous proofs, we will recover the Fourier coefficients of each
$A_j$, and in that manner, verify that $\psi$ is faithful.
Following the proof of Theorem~\ref{upper}, it suffices to show that
the functions
\[
 f_{v_s,\lambda}^k(A_j) = \ip{ \psi_{v_s,\lambda}(A_j)^* h_1,h_k }
\]
can be recovered from $\psi(A)$ for all $j \in J$ and $k\ge1$.

Fix $j_0 \in J$.
Since $A_{j_0} = P_{y_{j_0}} A_{j_0}$,
$f_{v_s,\lambda}^k(A_{j_0}) = 0$ if $r(v_s) \ne y_{j_0}$.
If $r(v_s) = y_{j_0}$ and $\lambda \in \bT^{|v|+1}$ is given, then there
is a sequence $i_n \in V_{j_0}$ so that $u_{i_n} = v_{s,k}$ and
$\lim_{n\to\infty} \lambda_{v_s,k} = \lambda$.
The result will follow if we establish that
\[
 \ol{ f_{v_s,\lambda}^k(A_{j_0}) } =
 \lim_{n\to\infty} 2^{i_n} \ip{ \psi(A)^* h_i, h_{s_{i_n}} } .
\]

We have observed that $\psi(A)^* h_i = \sum_{j \in J} \psi(A_j)^* k_j$.
So
{\allowdisplaybreaks
\begin{align*}
  2^{i_n} \ip{ \psi(A)^* h_i, h_{s_{i_n}} } &=
  2^{i_n} \sum_{j \in J} \ip{ \psi(A_j)^* k_j, h_{s_{i_n}} } \\
  &= 2^{i_n} \sum_{j \in J} \ip{ k_j, \psi(A_j) h_{s_{i_n}} } \\
  &= \ip{ h_{r_{i_n}}, \psi(A_j) h_{s_{i_n}} } +
     \sum_{j \in J} \sum_{\substack{1 \le i < i_n \\ i \in V_j}}\!\!
     2^{i_n - i} \ip{ h_{r_i}, \psi(A_j) h_{s_{i_n}} }
\end{align*}
} 
The first term has the desired limit.
So it remains to show that the second sum tends to zero.
The key is an estimate of $|\ip{ h_{r_i}, \psi(A_j) h_{s_{i_n}} }|$.

Observe that we can factor the representation $\psi_{w,\lambda}$ of $\fL_{G_0}$
as $\sigma \rho$ where $\rho$ replaces the maps
\[
 \psi_{w,\lambda}(e) = \frac 12 \sum_{e_{j_k}=e}h_{k-1} h_k^*
 \qfor e \in \E(G_0) \bsl \L
\]
with
\[
 \rho(e) = \frac 1 {\sqrt 2} \sum_{e_{j_k}=e}h_{k-1} h_k^*
 \qfor e \in \E(G_0) \bsl \L .
\]
and leaves
\[
 \rho(x) = \psi_{w,\lambda}(x) \qand
 \rho(f_x) = \psi_{w,\lambda}(f_x) \qfor x \in \V(G_0) .
\]
Then we define $\sigma$ to be the endomorphism of the nest algebra into
itself by multiplying the $n$th diagonal entries by $2^{-n/2}$.
Then $\rho$ is completely contractive by Corollary~\ref{dilationcor}.

We can now estimate $|\ip{ h_{r_i}, \psi(A_j) h_{s_{i_n}} }|$
by observing that this is a matrix entry of $\psi_{w,\lambda}(A_j)$
lying on the $s_{i_n} - r_i$ diagonal, and it is therefore bounded by
$\|A_j\| 2^{(r_i - s_{i_n})/2}$.
Since $\|A_j\| \le \|A\|$ and $r_{i-1} - r_i \ge 4i$,
\begin{align*}
  \sum_{j \in J} \sum_{\substack{1 \le i < i_n \\ i \in V_j}}
     2^{i_n - i} |\ip{ h_{r_i}, \psi(A_j) h_{s_{i_n}} }|
  &\le \|A\| \sum_{i=1}^{i_n-1} 2^{i_n - i} 2^{-\sum_{i<k\le i_n} 2i} \\
  &\le 2^{-i_n} \|A\| .
\end{align*}
The desired limit is now evident.
\bx

\begin{rem}\label{Zrepn}
It is clear that taking the transpose of the graphs yields a characterization
of which nests have a faithful representation onto the lower triangular nest
algebra.  Moreover it is not too difficult to formulate the corresponding
result for maximal nests of order type $\bZ$.  In particular, a strongly
transitive graph with a loop edge at every vertex has such a representation.
\end{rem}

As we can see, the representations constructed in this section are not
$*$-extendible.
This is not a limitation of our construction when the nest is non-trivial.
Indeed if $\fA_k$ is the norm closed algebra of the graph $P_k$ with one
vertex and $k$ loop edges, the non-commutative disk algebra of order
$k$, we have the following result.

\begin{prop} \label{P:not_*-extendible}
Let $\E_k$ be the C*-algebra generated by $\lambda(\fA_k)$, $k\geq2$,
and let $\phi$ be a $*$-representation of $\E_k$ in $\B(\H )$ that maps
$\fA_k$ faithfully into a nest algebra $\Alg \N$, with $\N \ne \{0,I\}$.
Then $\phi (\fA_k )$ is not dense in $\Alg \N$.
\end{prop}

\Prf Note that $\phi$ maps the creation operators onto isometries with
orthogonal ranges.  Therefore, the weakly closed algebra
$\fS = \overline{\phi(\fA_k)}$ generated by $\phi (\fA_k )$ is a free
semigroup algebra.  Hence, \cite[Corollary 2.9]{DKP} shows that the
Jacobson radical of $\fS$ consists of nilpotent operators of order at
most two.

The only non-trivial nests with this property have exactly three
elements, $\N = \{ 0, N, \H \}$.   From the Structure Theorem for free
semigroup algebras \cite{DKP}, this would imply that $P(N)^\perp$ is the
structure projection $P$. But $P(N)$ is another projection with the
property $P(N) \T(\N) P(N)$ is self-adjoint.  For free semigroup
algebras, $P$ is the unique maximal projection with this property. So a
free semigroup algebra cannot be a nest algebra other than $\B(\H)$.
\bx


\end{document}